\documentclass[12pt,a4paper]{article}
\usepackage{color}
\usepackage{amsmath}
\usepackage{amsfonts}
\usepackage{amssymb}
\usepackage{latexsym}
\usepackage{latexsym}
\usepackage{array}

\newtheorem{Theorem}{Theorem}

\newtheorem{Prop}{Proposition}

\newcommand{\qed}{\ifhmode\unskip\nobreak\fi\hfill \ifmmode\Box\else$\mathsurround=0pt\Box$\fi}

\voffset=-.5in \hoffset=-.5in \numberwithin{equation}{section}

\begin{document}
\title {
 Homogenization of Fractional Kinetic Systems
with Random Initial Data}

%\vskip 18 pt

\author { Gi-Ren Liu and Narn-Rueih Shieh\thanks{Correspondent author; E-mail: shiehnr@math.ntu.edu.tw}
\thanks{ Research partially supported by
a Taiwan NSC grant.}
\\
Department of Mathematics, National Taiwan University
\\
Taipei 10617, Taiwan
\\
% E-mail: shiehnr@math.ntu.edu.tw
}
 \maketitle
%\vskip 18 pt

 %{\it Short title:} Scaling limits for diffusion-wave systems

%\vskip 18 pt

\begin{abstract}
{\small Let ${\bf{w}}(t,x):=(u,v)(t,x),\ t>0,\ x\in
\mathbb{R}^{n},$ be the $\mathbb{R}^2$-valued spatial-temporal
random field ${\bf{w}}=(u, v)$ arising from a certain two-equation
system of fractional kinetic equations of reaction-diffusion type,
with  given random initial data $u(0,x)$ and $v(0,x).$ The
space-fractional derivative is characterized by the composition of
the inverses of the Riesz potential and the Bessel potential. We
discuss two scaling limits, the macro and the micro, for the
homogenization of ${\bf{w}}(t,x)$, and prove that the rescaled
limit is a singular field of multiple It\^{o}-Wiener   integral
type, subject to suitable assumptions  on the random initial
conditions. In the two scaling procedures, the Riesz and the
Bessel parameters play distinctive roles. Moreover, since the
component fields $u,v$ are dependent on the interactions present
within the system, we employ a certain stochastic decoupling
method to tackle this components dependence.  The time-fractional
system  is also considered, in which the Mittag-Leffler function
is used.

%The work shows, in particular, the various non-Gaussian scenarios
%proposed by works of Leonenko {\it et.al.}
%\cite{AL,AnhKinetic,AnhHomo,AnhSpectral,KLR,LM,LW} and the
%references therein, for the single diffusion type equations,  in
%classical or in fractional time/space derivatives, can be
%discussed for the two-equation system, in a significant way.
%Indeed, our result on the micro-scaling is new in the single
%equation  case.
}

\end{abstract}

\vskip 18 pt

{\it 2000 AMS classification numbers:} 60G60; 60H05; 62M15;
      35K45.


\bigskip

%%
%%%%%%%%%%%%%%%%%%%%%%%%%%%%%%%%%%%%%%%%%%%%%%%%%%%%%%%%%%%%%%%%%%%%%%%%%%%%%%%%%%%%%%%%%%%%%
%%%%%%%%%%%%%%%%%%%%%%%%%%%%%%%%%%% 1.Introduction %%%%%%%%%%%%%%%%%%%%%%%%%%%%%%%%%%%%%%%%%
%%%%%%%%%%%%%%%%%%%%%%%%%%%%%%%%%%%%%%%%%%%%%%%%%%%%%%%%%%%%%%%%%%%%%%%%%%%%%%%%%%%%%%%%%%%
\section{Introduction}
The purpose of  this paper is to  present a certain homogenization
theory associated with the following  linear fractional kinetic
system of reaction-diffusion type
\begin{align}\label{generalsystem}
\left\{
\begin{array}{lll}
\frac{\partial^{\beta}}{\partial
t^{\beta}}u(t,x)=-\mu_{1}(I-\Delta)^{\frac{\gamma}{2}}
(-\Delta)^{\frac{\alpha}{2}} u(t,x)+b_{11}u+b_{12}v,
& u(0,x)=u_{0}(x), \\
\frac{\partial^{\beta}}{\partial
t^{\beta}}v(t,x)=-\mu_{2}(I-\Delta)^{\frac{\gamma}{2}}
(-\Delta)^{\frac{\alpha}{2}} v(t,x)+b_{21}u+b_{22}v , &
v(0,x)=v_{0}(x),
\end{array}
\right.
\end{align}
 in the above, $\mu_i>0,\, 0< \beta \le 1, 0<\alpha\le 2, 0\le \gamma,$ and $ t>0,x\in
\mathbb{R}^{n}$. The parameter $\beta$ denotes the time-fractional
index,  and $\alpha, \gamma$ denote the space-fractional indices,
for which we refer as the Riesz parameter and the Bessel parameter
respectively (see \cite[V.1 and V.3]{Stein}).

When $(\beta,\alpha,\gamma)=(1,2,0)$, the system
(\ref{generalsystem}) is reduced to a classical reaction-diffusion
system. The time-fractional index $\beta< 1$ means sub-diffusive
(super-diffusive in case $\beta>1$, which we do not study in this
paper; see the remark in Section 6). The spatial-fractional Riesz
index $\alpha$ means the jumps of the evolution, and Bessel index
$\gamma$ means the tempering of large jumps; see the now-classic
book of Stein \cite[V.1 and V.3]{Stein} for precise mathematical
explanations. To our knowledge, fractional kinetic equations of
Riesz-Bessel type appear firstly in Anh and Leonenko
\cite{AnhKinetic,AnhHomo}; abundant subsequent works in this
direction by the authors and collaborators can be seen in
\cite{AnhSpectral,AnhHarmonic,Higherorder,BL,KLR,LR} and the
references therein. The two papers \cite{BL,LR} with external
potentials are particularly related to the study of this paper. We
should mention that fractional operators with two fractional
parameters are natural mathematical objects to describe long-range
dependence and/or intermittency; one can find data exhibiting such
characteristics in a large number of fields including economics,
finance, telecommunications, turbulence, and hydrology.

In this paper, we consider the system (\ref{generalsystem}) with
$\mu_1=\mu_2=\mu
>0$ and  with the random initial data $u_{0}$ and $\ v_{0}$,
of which are independent and each one has a certain long-range
dependence in its random structure. This paper is  along the
\cite{AnhHomo} on single fractional kinetic equation; yet our
results in this paper are with the following novel features.
Firstly and most importantly, we study two scaling procedures, the
macro and the micro, of the homogenization of the associated
spatial-temporal random solution-field, in which  the Riesz and
the Bessel parameters play distinctive roles; our result on the
micro-scaling is new, even for the single equation case; this
micro-scaling makes use of both the Riesz and the Bessel indices
and also needs the rescaling on the initial data. We feel that our
result in this micro-scaling may capture the proclaimed feature of
the intermittency of the random motions. Moreover, due to the
interactions present within the system, the components fields
$u,v$ are dependent (even we have assumed the independence of the
initial data $u_0,v_0$), and we employ a certain stochastic
decoupling method to tackle this components dependence. Our study
may show how the theory of Riesz potentials and Bessel potentials,
as in the Chapter V of Stein \cite{Stein}, may appear
significantly in the homogenization of random fields.

The study on the single P.D.E. with random initial condition can
be traced back to \cite{KF} and \cite{R}, and then has a long
active development; we refer to the citations in the above and the
references therein. There also has very significant progress on
Burgers' equation with random initial data; see the monograph of
Woyczy\'{n}ski \cite{W lecture} and the Chapter 6 of Bertoin
\cite{B}. Whilst, to our knowledge, relevant study on P.D.E.
system with random initial data seems few in previous literatures,
except the works of Leonenko and Woyczy\'{n}ski \cite{LW1,LW2} on
multi-dimensional Burgers' random fields (Burgers' Turbulence).

The results of this paper show that, for  the fractional kinetic
system (\ref{generalsystem}) with suitable random initial data,
the rescaled field ${\bf{w}}^{\varepsilon}(t,x)$ in the
micro-scaling (which we mean  $\varepsilon t$ and
$\varepsilon\downarrow 0$) of homogenization both the Riesz
parameter $\alpha $ and the Bessel parameter $\gamma $ play their
roles; while in the macro-scaling (which we mean
$\frac{t}{\varepsilon}$ and $\varepsilon\downarrow 0$) only the
Riesz parameter $\alpha $ plays the role. Nevertheless, in either
case the limiting field is a singular field of multiple
It\^{o}-Wiener integral type. Furthermore, the component fields
$u,v$ are dependent, due to the interactions present within the
system, we employ a certain {\it stochastic decoupling} method to
tackle this components dependence. The method  itself could be
potentially important in the future study on some random systems,
for example the gradient system of Hamilton-Jacobi equation with
random initial data, as in \cite[p.173]{W lecture} ; we notice
that the decoupling has been traditionally used in solving
differential equation systems.

The underlying idea in this paper is motivated by those works in
\cite{AnhKinetic,AnhHomo,BL,LR,LW} and the references therein.
Namely, we use
 the spectral representations to describe the sample field arising from the initial data,
and the relations between Hermite polynomials and homogeneous
chaos associated with the initial data, to get representations for
the limit field in terms of multiple It\^{o}-Wiener integrals.
From limit theorems point-of-view, our results, also those in the
above citations, belong to the realm of non-central limit theorems
for convolution type integrals, in which the papers \cite{Taq, DM}
are pioneering; see  also the monograph of Major\cite{M} and
survey papers in the special volume edited by Doukhan, Oppenheim
and Taqqu \cite{DOT}

The paper is organized as follows. In Section 2 we give the
explicit solution of the system (\ref{generalsystem}). In Section
3 the initial data are assumed to be stationary random fields, and
we discuss the covariance structure of the resulting
solution-vector random field of (\ref{generalsystem}), subject to
the specified random initial condition; we show that the spectral
method is suitable in describing our random field relative to the
space-time parameter. We also introduce the initial field to be a
certain subordinated Gaussian random field generated by
 a class of non-random functions whose variables are relative to the spatial parameter.
In the main Sections 4 and 5, we consider (\ref{generalsystem})
with the usual time-derivative and the fractional
spatial-derivative characterized by the Riesz and the Bessel
parameters. We present the homogenization of micro-scaling  in
Section 4 and the less subtle macro-scaling in section 5,
respectively. In Section 6 we provide  extensions of the results
in Sections 4 and 5 to the time-fractional $\beta < 1$, in which
we need to use the Mittag-Leffer function. The proofs of all our
results are given in Section 7.

{\bf Acknowledgement.} The authors are grateful to the inspiring
lectures of Professor W.A. Woyczy\'{n}ski at National Taiwan
University for the perspective on {\it Mathematical Theory} of
Fractional P.D.E.

%{\bf Remark on the context:} in the followings we mainly discuss the
%range $\beta\in (1,2)$.
% though some statements in the
%preliminary Section 2 is written for the full range $\beta\in
%(0,2]$ . This is based on the following reason.
%The case $\beta =1$, that is the classical reaction-diffusion case,
%has been treated in \cite{LSh08}. The case $\beta=2$, that is the
%classical wave type system, must be treated in the framework of
%generalized random fields, and thus we treat it in a separate space.
%The case $\beta\in (0,1)$ can be essentially treated in the same way
%as $1<\beta<2$, and is simpler since there are no initial data
%$u_t,v_t$ involved; thus we leave it to the reader, see the remark
%below Theorem 2.
%In the physics language,
%for $\beta\in (1,2)$ we are dealing with the sub-oscillating
%diffusion and hyperbolic time-fractional reaction-diffusion
%system, see the remark in \cite{GDM} at the end of Section 1
%(p.216) there.
%Moreover, we choose the spatial dimensionality to be 3, which may
%properly cover all the various possible cases mentioned in Theorem
%1; see the remark below Theorem 1. While other dimensionality, say
%$x\in \mathbf{R}^{n }, n=1,2$, does not exactly, though does
%partially, meet this aspect; thus we also leave  it to the reader.
%%%%%%%%%%%%%%%%%%%%%%%%%%%%%%%%%%%%%%%%%%%%%
\vskip 20 pt

\section{Preliminaries}
To begin with, we rewrite the system (\ref{generalsystem}), with
 $\mu_1=\mu_2=\mu
>0$, $\alpha>0,\ \gamma\geq 0$ and $\beta=1$,
in the matrix form as follows:
\begin{align}\label{system}
 \frac{\partial}{\partial t} \left(\begin{array}{cccc}u\\v\end{array}\right)=
 -\mu(I-\Delta)^{\frac{\gamma}{2}}
(-\Delta)^{\frac{\alpha}{2}}\left(\begin{array}{cccc}u\\v\end{array}\right)+
B\left(\begin{array}{cccc}u\\v\end{array}\right),
\end{align}
subject to some initial conditions
\begin{align}\label{initial}
\left(\begin{array}{cccc}u(0,x)\\v(0,x)\end{array}\right)=\left(\begin{array}{cccc}u_{0}(x)\\v_{0}(x)\end{array}\right),
\ x\in \mathbb{R}^{n},
\end{align}
where $u=u(t,x)$, $v=v(t,x),\ t>0$, $x\in \mathbb{R}^{n}$, $\Delta$
is the $n$-dimensional Laplacian, and $B$ is a $2\times2$ matrix.

The Green function $G(t,x;\alpha,\gamma)$ associated with the
operator $\partial_t + \mu(I-\Delta)^{\frac{\gamma}{2}}
(-\Delta)^{\frac{\alpha}{2}}$ is represented via the
  spatial Fourier
transform as follows; see, \cite[Chapter 5]{Stein} or
\cite[Section 2]{AnhHomo}.
\begin{align}\label{fractional Green function}
\int_{\mathbb{R}^{n}}e^{i<x,\lambda>}G(t,x;\alpha,\gamma)dx
=\textup{exp}[-\mu t |\lambda|^{\alpha}(1+|\lambda|^{2})^{\frac{\gamma}{2}}],\
\lambda\in \mathbb{R}^{n},
\end{align}
where
$<\cdot,\cdot>$ denotes the inner product on
$\mathbb{R}^{n}$.

In order to get a explicit representation for the solution of (\ref{system}),
we impose the following assumption on the matrix $B$.\\
%%%%%%%%%%%%%%%%%%%%%%%%%%%%%%%%%%%%%%%% Condition A. %%%%%%%%%%%%%%%%%%%%%%%%%%%%%%%%%
\begin{bfseries}Condition A.\end{bfseries}  Suppose the
matrix $[b_{ij}]_{1\leq i,j\leq 2}$ is  diagonalizable, i.e., the
matrix $B$ can be written as
\begin{align}\label{matrixB}
B=\left(\begin{array}{cc} b_{11} & b_{12} \\b_{21} & b_{22}
\end{array} \right)=PDP^{-1},\ \textup{with}\
P=\left(\begin{array}{cc}
 p_{11} & p_{12} \\p_{21} & p_{22} \end{array} \right),
\end{align}
where $P$ is a real-valued non-degenerate eigenvector matrix
associated with the matrix $B$, and $D=diag(d_{1},d_{2})$,
$d_{1},\ d_{2}\in \mathbb{R}$, where $d_j$ is the eigenvalue
associated with the eigenvector $(p_{1,j}, p_{2,j})^T$ (here and
henceforth, T denotes the transpose). Without loss of generality,
we suppose that $det(P)=1$.

Under Condition A, the Cauchy problem (\ref{system})
(\ref{initial}) has the unique solution given by
\begin{align}\label{solution}
\left(\begin{array}{cccc}u(t,x)\\v(t,x)\end{array}\right)=
Q(t;d_{1},d_{2})\left(\begin{array}{cccc}U(t,x)\\V(t,x)\end{array}\right),\
 t>0,\ x\in \mathbb{R}^{n},
\end{align}
where
\begin{align}\label{Q matrix}
Q(t;d_{1},d_{2}):=P\left(\begin{array}{cccc}e^{d_{1}t} & 0 \\0 &
e^{d_{2}t} \end{array}\right)P^{-1},
\end{align}
and $U(t,x)$, $V(t,x)$  are determined by
\begin{align}\label{homogeneous case}
\left(\begin{array}{cccc}U(t,x)\\V(t,x)\end{array}\right)=\int_{\mathbb{R}^{n}}G(t,y;\alpha,\gamma)\left(\begin{array}{cccc}
u_{0}(x-y)\\v_{0}(x-y)\end{array}\right)dy,
\end{align}
where the Green function $G(t,y;\alpha,\gamma)$ is defined in
(\ref{fractional Green function}).
%by the
%transformation
%\begin{align}\label{Green}
%G(t,y;\alpha,\gamma)=(2\pi)^{-n}\int_{\mathbb{R}^{n}}e^{-i<\lambda,y>}e^{-\mu
%t|\lambda|^{\alpha}(1+|\lambda|^{2})^{\frac{\gamma}{2}}}d\lambda.
%\end{align}

For completeness, we give the sketchy proofs of (\ref{solution}).
Firstly, by taking the spatial Fourier transform on both sides of
(\ref{system}), under Condition A, we have
\begin{align*}
\frac{\partial}{\partial t}
&\left(\begin{array}{cc}\widehat{u}\\\widehat{v}\end{array}\right)(t,\lambda)
=(-\mu|\lambda|^{\alpha}(1+|\lambda|^{2})^{\frac{\gamma}{2}}
+B)
\left(\begin{array}{cc}\widehat{u}\\\widehat{v}\end{array}\right)(t,\lambda)\\
&= P
\begin{pmatrix}-\mu |\lambda|^{\alpha}(1+|\lambda|^{2})^{\frac{\gamma}{2}}+d_{1}
& 0
\\0 & -\mu |\lambda|^{\alpha}(1+|\lambda|^{2})^{\frac{\gamma}{2}}+d_{2}\end{pmatrix}
P^{-1}
\left(\begin{array}{cc}\widehat{u}\\\widehat{v}\end{array}\right)(t,\lambda)
\end{align*}
Thus,
\begin{align}\notag
&\left(\begin{array}{cc}\widehat{u}\\\widehat{v}\end{array}\right)
(t,\lambda)
\\\notag
&=\textup{exp}\Big\{t P
\begin{pmatrix}-\mu|\lambda|^{\alpha}(1+|\lambda|^{2})^{\frac{\gamma}{2}}+d_{1}
& 0
\\0 & -\mu|\lambda|^{\alpha}(1+|\lambda|^{2})^{\frac{\gamma}{2}}+d_{2}\end{pmatrix}
P^{-1} \Big\}
\left(\begin{array}{cc}\widehat{u}\\\widehat{v}\end{array}\right)(0,\lambda)\\\label{Appendix
1} &= P\left(\begin{array}{cc}e^{d_{1}t} & 0\\0 &
e^{d_{2}t}\end{array}\right)P^{-1} e^{-\mu
t|\lambda|^{\alpha}(1+|\lambda|^{2})^{\frac{\gamma}{2}}}
\left(\begin{array}{cc}\widehat{u}\\\widehat{v}\end{array}\right)(0,\lambda).
\end{align}
Finally, (\ref{solution}) and (\ref{homogeneous case}) are followed by taking
the inverse Fourier transform on both sides of (\ref{Appendix
1}) and using the representation (\ref{fractional Green function}).
Additionally, by (\ref{fractional Green function}) we can also observe that
\begin{align}\label{integral of G}
\int_{\mathbb{R}^{n}}G(t,x;\alpha,\gamma)\ dx = 1,\ \textup{for any}\ t\geq 0.
\end{align}

%%%%%%%%%%%%%%%%%%%%%%%%%%%%%%%   Section : Inter-relative random structure %%%%%%%%%%%%%%%%%%%%%%%%%
\section{Correlated random structures}

\subsection{general random initial data}
Firstly, we set $(\Omega, \mathcal{F}, \mathcal{P})$ to be an
underlying probability space,
 such that all random element appeared in this paper are measurable with respect to it.\\
The following condition is imposed on the initials, in which and
henceforth.\\
%%%%%%%%%%%%%%%%%%%%%%%%%%%%%%%%%%%%%%% Condition B. %%%%%%%%%%%%%%%%%%%%%%%%%%%%%%%%%
\begin{bfseries}Condition B.\end{bfseries}
Let $u_{0}(x)=\eta_{1}(x)=\eta_{1}(x,\omega)$ and
$v_{0}(x)=\eta_{2}(x)=\eta_{2}(x,\omega),\ x\in\mathbb{R}^{n},
\omega\in\Omega$, be independent, and distributed as two real,
mean-square continuous, homogeneous and isotropic random fields
defined on the underlying complete probability space
$(\Omega,\mathcal{F},\mathcal{P})$. We assume that $
E\eta_{j}(x)\equiv 0,\ \textup{Var}(\eta_{j}(x))\equiv 1, $ and covariance
functions
\[
R_{\eta_{j}}(x)=\widetilde{R}_{\eta_{j}}(|x|):=\textup{Cov}(\eta_{j}(0),\eta_{j}(x))=\int_{\mathbb{R}^{n}}e^{i<\lambda,x>}F_{j}(d\lambda),\
\ j\in \{1,2\},
\]
where the last equality is guaranteed by Bochner-Khintchine theorem
and $F_{j}(\cdot)$ is the spectral measure corresponding to the
field $\eta_{j}(\cdot)$ for $j\in\{1,2\},$ respectively.
\\
%%%%%%%%%%%%%%%%%%%%%%%%%%%%%%%%%%% End of  Condition B. %%%%%%%%%%%%%%%%%%%%%%%%%%%%%%%%%%%%%%
\indent Under Condition B, in view of Karhunen's Theorem (see, for
example,  Gihman and Skorokhod \cite{GS}, pp. 208-230), there exist
complex-valued orthogonally scattered random measures $Z_{F{j}},\
j\in\{1,2\}$, such that the random fields $\eta_{j}(x),\ j\in\{1,2\}$, have the following
spectral representations
\begin{align}
\eta_{j}(x)=\int_{\mathbb{R}^{n}}
e^{i<\lambda,x>}Z_{F_{j}}(d\lambda),\ \ \
j\in\{1,2\},\label{Spectral represent}
\end{align}
where $EZ_{F_{k}}(\Delta_{1})=0,\
EZ_{F_{k}}(\Delta_{1})\overline{Z_{F_{j}}}
(\Delta_{2})=\delta^{j}_{k}F_{j}(\Delta_{1}\cap\Delta_{2})$, for any
$j,\ k\in\{1,2\}$ and $\Delta_{1},\Delta_{2}\in
\mathcal{B}(\mathbb{R}^{n})$ ($\delta_{k}^{j}$ is the Kronecker
symbol).

By the above spectral representation for the initial data
(\ref{initial}), we can describe the vector-solution
$\{(u(t,x),v(t,x)), t>0, \ x\in \mathbb{R}^{n}\}$ by stochastic
integration :
%%%%%%%%%%%%%%%%%%%%%%%%%%%%%%%%%%%%%% Proposition 1. %%%%%%%%%%%%%%%%%%%%%%%%%%%%%%%%%%%%%%%%%%%%%
\begin{Prop}
Let
$\mathbf{w}(t,x;\mathbf{w}_{0}(\cdot)):=(u(t,x;u_{0}(\cdot),v(t,x;v_{0}(\cdot)),\
\mathbf{w}_{0}(\cdot)=(u_{0}(\cdot),v_{0}(\cdot))$, be the
vector-solution of the initial value problem (\ref{system})
(\ref{initial}), of which satisfies Condition A and B, then
\begin{align}\label{spectral solution}
\mathbf{w}(t,x;\mathbf{w}_{0}(\cdot)) =Q(t;d_{1},d_{2})
\int_{\mathbb{R}^{n}} e^{i<\lambda,x>} e^{-\mu
t|\lambda|^{\alpha}(1+|\lambda|^{2})^{\frac{\gamma}{2}}}
\left(\begin{array}{cccc}
Z_{F_{1}}(d\lambda)\\Z_{F_{2}}(d\lambda)\end{array}\right),
\end{align}
with the inter-relative covariance structure
\begin{align}\notag
E\mathbf{w}(t,x;\mathbf{w}_{0}(\cdot))\mathbf{w}^{T}(t^{'},x^{'};\mathbf{w}_{0}(\cdot))
\end{align}
\begin{align}\label{covariance structure}
=\int_{\mathbb{R}^{n}}e^{i<\lambda,x-x^{'}>} e^{-\mu
(t+t^{'})|\lambda|^{\alpha}(1+|\lambda|^{2})^{\frac{\gamma}{2}}}
Q(t;d_{1},d_{2})
 \begin{pmatrix}F_{1}(d\lambda) & 0\\0 &
F_{2}(d\lambda)
\end{pmatrix}
Q(t^{'};d_{1},d_{2})^{T},
\end{align}\\[-0.7cm]
{\it where $Q(t;d_{1},d_{2})$ is defined in (\ref{Q matrix}).}
\end{Prop}
%%%%%%%%%%%%%%%%%%%%%%%%%%%%%%%%%%%%%%%%% End of Proposition 1. %%%%%%%%%%%%%%%%%%%%%%%%%%%%%%%%%%%

%%%%%%%%%%%%%%%%%%%%%%%%%%%%% 3.2 Subordinated Gaussian random field %%%%%%%%%%%%%%%%%%%%%%%%%%%%%%
%%%%%%%%%%%%%%%%%%%%%%%%%%%%%%%%%%%%%%%%%%%%%%%%%%%%%%%%%%%%%%%%%%%%%%%%%%%%%%%%%%%%%%%%%%%%%%%%
\subsection{Subordinated Gaussian initial data}
In this subsection, we assume further that the initials are subordinated fields, as follows:\\
%%%%%%%%%%%%%%%%%%%%%%%%%%%%%%%%%%%%%%%%% Condition C. %%%%%%%%%%%%%%%%%%%%%%%%%%%%%%%%%%%%%%%%%%%%%
\begin{bfseries}Condition C.\end{bfseries} We consider the random initial data (\ref{initial})
$\mathbf{w}_{0}(x):=(u_{0}(x),v_{0}(x))=(\eta_{1}(x),\eta_{2}(x)),\
x\in \mathbb{R}^{n},$ satisfies Condition B and each component has
the following form
\begin{align}\label{initial data form}
\eta_{j}(x):=h_{j}(\zeta_{j}(x)),\ x\in \mathbb{R}^{n},\ j\in
\{1,2\}.
\end{align}
The $\zeta_{1}(x)$ and $\zeta_{2}(x)$ are independent, mean-square
continuous, homogeneous and isotropic Gaussian random fields, each
is of mean zero and of variance 1, and for each the spectral
measure $F_{j}(d\lambda)$ has the (spectral) density
$f_{j}(\lambda),\ \lambda\in \mathbb{R}^{n}$, and $f_{j}(\lambda)$
is decreasing for $|\lambda|>\lambda_{0}$ for some $\lambda_{0}>0$
and continuous for all $\lambda\neq 0$, $j\in\{1,2\}$,
respectively. Moreover, we assume that $h_{j}(\cdot),\
j\in\{1,2\}$, are real non-random Borel functions satisfy
\begin{align}
Eh^{2}_{j}(\zeta_{j}(0))<\infty,\ \ j\in\{1,2\}.\label{L2 function}
\end{align}
%%%%%%%%%%%%%%%%%%%%%%%%%%%%%%%%%%%%%%%%%% End of Condition C. %%%%%%%%%%%%%%%%%%%%%%%%%%%%%%%%%%%%%%%%

Under Condition C, we have the spectral representations for the
sample paths of $\zeta_{j}(x),\ j\in\{1,2\},$ as below:
\begin{align}
\zeta_{j}(x)=\int_{\mathbb{R}^{n}}e^{i<x,\lambda>}\sqrt{f_{j}(\lambda)}W_{j}(d\lambda),\
x\in \mathbb{R}^{n},\ j\in\{1,2\},\label{sample path represent}
\end{align}
where $W_{j}(A)$ is a Gaussian noise measure,
% which can be expressed
%formally as
%\begin{align}
%W_{j}(A)=\int_{\mathbb{R}^{n}}\mathcal{F}^{-1}(\frac{1_{A}(\cdot)}
%{\sqrt{f_{j}(\cdot)}})(x)\zeta_{j}(x)dx,\ \textup{for any}\ A\in
%\mathcal{B}(\mathbb{R}^{n}).
%\end{align}
%We also remark that,
$W_{j}(A),\ A\in \mathcal{B}(\mathbb{R}^{n})$,
are centered Gaussian with
$EW_{i}(d\lambda)\overline{W_{j}(d\mu)}=\delta^{j}_{i}\delta(\lambda-\mu)d\lambda\
d\mu$.
% In the above and henceforth, $\delta(\cdot)$ denotes the
%Delta distribution, and $\mathcal{F}^{-1}$ denotes the inverse
%Fourier transform.
\\
%%%%%%%%%%%%%%%%%%%%%%%%%%%% About  Hermite polynomial ... %%%%%%%%%%%%%%%%%%%%%%%%%%%%%%%%%%%%%%%%
\noindent Due to (\ref{L2 function}) in Condition C, we can consider
the following orthogonal expansions of $h_{j}(u)$ in the Hilbert
space $L^{2}(\mathbb{R},p(u)du)$ with
$p(u)=\frac{1}{\sqrt{2\pi}}e^{-\frac{u^{2}}{2}}$:
\begin{align}
h_{j}(u)=C^{(j)}_{0}+\sum_{\sigma=1}^{\infty}{C^{(j)}_{\sigma}
\frac{H_{\sigma}(u)}{\sqrt{\sigma !}}}, \ \ j\in\{1,2\},
\end{align}
where
\begin{align}\label{hermitecoeff}
C^{(j)}_{\sigma}=\int_{\mathbb{R}}h_{j}(u)\frac{H_{\sigma}(u)}{\sqrt{\sigma!}}p(u)du,\
\ j\in\{1,2\},
\end{align}
and $\{H_{\sigma}(u),\sigma=0,1,2,...\}$ are the Hermite
polynomials, that is,
\[H_{\sigma}(u)=(-1)^{\sigma}e^{\frac{u^{2}}{2}}\frac{d^{\sigma}}{du^{\sigma}}e^{\frac{-u^{2}}{2}},\ \ \textup{for}\ \ \sigma\in \{0,1,2,...\}.\]
%%%%%%%%%%%
It is known that the following two important properties hold ( see,
for example, Major \cite{M}, Corollary 5.5 and p. 30 ) :
\begin{align}\label{expectionhermite}
E[H_{\sigma_{1}}(\zeta_{j}(y_{1}))H_{\sigma_{2}}(\zeta_{j^{'}}(y_{2}))]=
\delta^{j^{'}}_{j}\delta^{\sigma_{1}}_{\sigma_{2}} \sigma_{1}!
R^{\sigma_{1}}_{\zeta_{j}}(y_{1}-y_{2}),\ \ \ y_{1},y_{2}\in
\mathbb{R}^{n}
\end{align}
and
\begin{align}\label{itoformula}
H_{\rho}(\zeta_{j}(x))=
\int^{'}_{\mathbb{R}^{n\times\rho}}e^{i<x,\lambda_{1}+...+\lambda_{\rho}>}\prod_{\sigma=1}^{\rho}{\sqrt{f_{j}(\lambda_{\sigma})}}
W_{j}(d\lambda_{\sigma}).
\end{align}
In  the integral representation (\ref{itoformula}),  the
integration $\int^{'}$ means that it excludes the diagonal
hyperplanes $z_{i}=\mp z_{j}, i, j=1,..., \rho, i\neq j$.

%%%%%%%%%%%%%%%%%%%%%%%%%
\noindent The {\it Hermite rank } of the functions $h_{j}(\cdot)$ is
defined by
\begin{align*}
m_{j}:=\textup{inf}\{\sigma\geq 1:\ C^{(j)}_{\sigma}\neq 0\} ,\ \
j\in\{1,2\}.
\end{align*}
Specializing Proposition 1 in Subsection 3.1 to the present
subordinated Gaussian initials, we have
%%%%%%%%%%%%%%%%%%%%%%%%%%%%%%%%%%%%%%%%%%%%%%%%%%%%%%%%%%%%%%%%%%%%%%%%%%%%%%%%%%%
%%%%%%%%%%%%%%%%%%%%%%%%%%%%%%% Proposition 2 %%%%%%%%%%%%%%%%%%%%%%%%%%%%%%%%%%
\begin{Prop}\label{proposition 1''}
 {\it Let $\mathbf{w}(t,x;\mathbf{w}_{0}(\cdot)):=(u(t,x;u_{0}(\cdot)),v(t,x;v_{0}(\cdot))),\ t>0\,\ x\in
\mathbb{R}^{n}\}$  be the vector-solution (\ref{solution}) of the
initial value problem (\ref{system}) (\ref{initial}), of which
satisfies Condition A and C, then} {\it all the statements in
Proposition 1 remain valid, with (\ref{spectral solution}) is
expressed as}
\begin{align*}
{\bf{w}}(t,x;\mathbf{w}_{0}(\cdot))=&Q(t;d_{1},d_{2})\Big\{
\left(\begin{array}{cc}C^{(1)}_{0}\\C^{(2)}_{0}\end{array}\right)+
\\
&\underset{\rho\in \mathbb{N}}{\sum}
\int^{'}_{\mathbb{R}^{n\times
\rho}}
e^{i<x,\lambda_{1}+\cdots +\lambda_{\rho}> -\mu t
|\lambda_{1}+\cdots +\lambda_{\rho}|^{\alpha}(1+|\lambda_{1}+\cdots
+\lambda_{\rho}|^{2})^{\frac{\gamma}{2}}}
\left(\begin{array}{cc}Z^{(\rho)}_{F_{1}}(d\lambda)\\Z^{(\rho)}_{F_{2}}(d\lambda)\end{array}\right)\Big\},
\end{align*}
{\it where}
\begin{align*}
\left(\begin{array}{cc}Z^{(\rho)}_{F_{1}}(d\lambda)\\Z^{(\rho)}_{F_{2}}(d\lambda)\end{array}\right)
:=
\begin{pmatrix}
\frac{C^{(1)}_{\rho}}{\sqrt{\rho!}}
\overset{\rho}{\underset{\sigma=1}{\prod}}{\sqrt{f_{1}(\lambda_{\sigma})}}
W_{1}(d\lambda_{\sigma})
\\
\frac{C^{(2)}_{\rho}}{\sqrt{\rho!}}
\overset{\rho}{\underset{\sigma=1}{\prod}}{\sqrt{f_{2}(\lambda_{\sigma})}}
W_{2}(d\lambda_{\sigma})
\end{pmatrix},
\end{align*}
{\it and the coefficient} $C^{(j)}_{\rho},\ j\in\{1,2\}$ {\it is
defined in (\ref{hermitecoeff}).}
\end{Prop}
%%%%%%%%%%%%%%%%%%%%%%%%%%%%%% end of Proposition 2 %%%%%%%%%%%%%%%%%%%%%%%%%%%%%

\indent We also impose the following assumption which is related
to the long-range dependence
 of the underlying Gaussian fields $\zeta_{j}(x),\ j\in\{1,2\}$;
we refer to \cite{AH,DOT} for the notion and the literatures of
long-range dependence. In the following and henceforth, the
notation $ f(\cdot)\sim g(\cdot)$ means that the ratio
$f(\cdot)/g(\cdot)$ tends to 1, as the indicated variable
``$\cdot$" tends to $\infty$ or tends to 0, according to the
context.\\
%%%%%%%%%%%%%%%%%%%% Condition D. %%%%%%%%%%%%%%%%%%%%%%%%%%%%%%%%%%%%%%%%%%%%%%%%%%%
\begin{bfseries}Condition D.\end{bfseries}  The Gaussian
random fields $\zeta_{j}(x),\ j\in\{1,2\}$, in Condition C, have
their covariance functions to be regular varying at infinity in the
sense that:
\begin{align}\label{long}
R_{\zeta_{j}}(x)\sim \frac{L(|x|)}{|x|^{\kappa_{j}}},\ \
\textup{as}\ |x|\rightarrow \infty,\ 0<\kappa_{j}<\frac{n}{m_{j}},\
j\in\{1,2\},
\end{align}
where $L:(0,\infty)\rightarrow(0,\infty)$ is a slowly varying
function at infinity and is bounded on each finite interval;
recall that $L$ is said to be slowly varying at infinity if
$\textup{lim}_{y\rightarrow \infty}[L(cy)/L(y)]=1$ uniformly for
any $c\in(a,b),\ 0<a<b<\infty$. Here, the relation $``\sim"$ in
(\ref{long}) means that
$\underset{|x|\rightarrow\infty}{\textup{lim}}|x|^{\kappa_{j}}R_{\zeta_{j}}(x)/L(|x|)=1$;
the similar notation will be used in this section and also in Section 7.\\
%%%%%%%%%%%%%%%%%%%%%%%%%%% End of Condition D. %%%%%%%%%%%%%%%%%%%%%%%%%%%%%%
%
Under Condition D, by  a Tauberian  theorem (see, for example, the
book of Leonenko \cite[p. 66]{L}), the spectral density functions
of the random fields $\zeta_{j}(x)$, $j\in\{1,2\}$, are regular
varying near the origin as follows:
\begin{align}\label{tauberian}
f_{j}(\lambda)\sim
K(n,\kappa_{j})|\lambda|^{\kappa_{j}-n}L(|\frac{1}{\lambda}|), \
\textup{as}\ \lambda\rightarrow 0,\ j\in \{1,2\},
\end{align}
where the Tauberian constant
$K(n,\kappa_{j})=\frac{\Gamma(\frac{n-\kappa_{j}}{2})}{2^{\kappa_{j}}\pi^{\frac{n}{2}}
\Gamma(\frac{\kappa_{j}}{2})}$.\\
%%%%%%%%%%%%%%%%%%%%%%%%%%%%%%%%%%%%%%%%%%%%%%%%%%%%%%%%%%%%%%%%%%%%%%%%%%%%%%%%%%%%%%%%%%%%%%%%%%%%%%%%%%%%%%%%%%
We note that, for each natural number $\rho \ge 2$, the power of
the covariance function $(R_{\zeta_{j}}(x))^{\rho}$ itself is
still the covariance function of some random field, for which
there exists the corresponding spectral density function
$(f_{j})^{*\rho}(\lambda)$. Indeed, the function
$(f_{j})^{*\rho}(\lambda), \lambda\in \mathbb{R}^{n}$, is the
$\rho-$th convolution of $f_{j}(\lambda)$ defined as
\begin{align}\label{convolutionspectral}
(f_{j})^{*\rho}(\lambda) =\int_{\mathbb{R}^{n\times (\rho-1)}}
f_{j}(\lambda-\lambda_{1}) f_{j}(\lambda_{1}-\lambda_{2})\cdots
f_{j}(\lambda_{\rho-2}-\lambda_{\rho-1})f_{j}(\lambda_{\rho-1})\overset{\rho-1}{\underset{l=1}{\prod}} d\lambda_{l}.
\end{align}
for $\rho\geq 2$. Since $L^{\rho}(|x|)$ is still a slowly varying
function for any $\rho$, when the
 $\rho$ satisfies $0<\rho\kappa_{j}<n$, we can apply
the Tauberian theorem again to get
\begin{align}\label{tauberian2}
(f_{j})^{*\rho}(\lambda)\sim
K(n,\rho\kappa_{j})|\lambda|^{\rho\kappa_{j}-n}L^{\rho}(|\frac{1}{\lambda}|),
\ \textup{as}\ \lambda\rightarrow 0,\ 0<\rho\kappa_{j}<n,
\end{align}
for $j\in \{1,2\}.$ While, if $\rho\kappa_{j}>n$ then the
covariance function $(R_{\zeta_{j}}(x))^{\rho}$ belongs to the
class $L^{1}(\mathbb{R}^{n})$; thus the corresponding spectral
density function is everywhere continuous and satisfies
\begin{align}\label{L1spectral}
(2\pi)^{n}(f_{j})^{*\rho}(0)=
\int_{\mathbb{R}^{n}}(R_{\zeta_{j}})^{\rho}(x)dx \leq
\int_{\mathbb{R}^{n}}|R_{\zeta_{j}}(x)|^{\rho}dx \leq
\int_{\mathbb{R}^{n}}|R_{\zeta_{j}}(x)|^{\rho^{*}}dx<\infty,
\end{align}
where $\rho^{*}:=\textup{inf}\{\rho\in \mathbb{N}|\
\rho\kappa_{j}>n\}$; we note that $|R_{\zeta_{j}}(\cdot)|\leq 1$.

The displays (\ref{convolutionspectral}), (\ref{tauberian2}) and
(\ref{L1spectral}) will be used  in the proofs in Section 7.

\vskip 20 pt
%%%%%%%%%%%%%%%%%%%%%%%%%%%%%%%%%%%%%%%%%%%%%%%%%%%%%%%%%%%%%%%%%%%%%%%%%%%%%%%%%%%%%%%%%%%%%%%%%%%%%%%%%%%%%%%%%%

%%%%%%%%%%%%%%%%%%%%%%%%%%%%%%%%%%% End of  Theorem 1. %%%%%%%%%%%%%%%%%%%%%%%%%%%%%%%%%%

\section{Micro-scalings for the solution vector-field }
In this section, we present the main result of this paper, which
concerns with the micro-scaling of the homogenization of the
spatial-temporal random field associated with (2.1), with the
initial data (2.2) subject to the conditions in Section 3. We show
that both the Riesz parameter $\alpha$ and the  Bessel parameter
$\gamma$ plays their roles in the scaling procedure. The results
in this section are more subtle than the macro-scaling discussed
in the next section; see the remark below Theorem 1 for the
interpretation.

Firstly, we prove the following micro-scaling of homogenization
for a single fractional kinetic equation, subject to the random
initial data.
\begin{equation}\label{single equation}
\frac{\partial s}{\partial
t}(t,x)=-\mu(I-\Delta)^{\frac{\gamma}{2}}
(-\Delta)^{\frac{\alpha}{2}} s(t,x),\ \ s(0,x)=h(\zeta(x)).
\end{equation}
To our knowledge,  the homogenization present in the below is a
completely new type. In (4.2) below, the notation imposed on
$\zeta$ wants to mean that the variable of $\zeta$ is under the
indicated dilation factor.
%%%%%%%%%%%%%%%%%%%%%%%%   prop. (micro) single case %%%%%%%%%%%%%%%%%%%%%%%%%%%%%%%%%%%%%%%%%%%%%%%%%%%%
\begin{Theorem}\label{Micro Scaling Theorem}
Let $s:=s(t,x;s_{0}(\cdot)),\ t> 0,\ x\in \mathbb{R}^{n},$ be a
solution of (\ref{single equation}), which satisfies the above Condition B, C
and D with $\kappa \in(0, \frac{n}{m})$,
where $m$ denotes the Hermite rank of
the non-random function $h(\cdot)$ on $\mathbb{R}$, which has the Hermite coefficients
$C_{i}(h),\
i=0,1,\ldots$ (i.e., $h_{1}(x)=h(x),\
\zeta_{1}(x)=\zeta(x),\ f_{1}(\lambda)=f(\lambda)$ and
$\kappa_{1}=\kappa$, etc. in Section 3).
Then, for any fixed parameter $\chi>0$,\\
(1) The behaviour of the covariance function of the rescaled random
field $s^{\varepsilon}(t,x),\ t>0,\ x\in \mathbb{R}^{n}$,
\begin{align}\label{rescaled micro field}
s^{\varepsilon}(t,x):=[\varepsilon^{m\kappa\chi}L^{m}(\varepsilon^{-\chi})]^{-\frac{1}{2}}
\Big\{s(\varepsilon t, \varepsilon^{\frac{1}{\alpha+\gamma}}x
;h(\zeta(\varepsilon^{-\frac{1}{\alpha+\gamma}-\chi}\cdot)))-C_{0}(h)
\Big\},
\end{align}
is given by:
\begin{align}\label{cov micro single thm}
\underset{\varepsilon\rightarrow
0}{\textup{lim}}\textup{Cov}(s^{\varepsilon}(t,x)s^{\varepsilon}(t^{'},x^{'}))
=
(C_{m}(h))^{2}K(n,m\kappa)\underset{{\mathbb{R}^{n}}}{\int} e^{i<x-x^{'},\tau>}
 \frac{e^{ -\mu
(t+t^{'})
|\tau|^{\alpha+\gamma}}}{|\tau|^{n-m\kappa}}d\tau.
\end{align}
\\
(2) When $\varepsilon\rightarrow 0$, the rescaled random field
$s^{\varepsilon}(t,x),\ t>0,\ x\in \mathbb{R}^{n},$ converges to the
limiting spatial-temporal random field $s_{m}(t,x),\ t>0,\ x\in
\mathbb{R}^{n}$, in the finite dimensional distribution sense, and
$s_{m}(t,x)$ is represented by the Multiple-Wiener integrals
\begin{align}\label{limiting field(micro)}
s_{m}(t,x):=\frac{C_{m}(h)}{\sqrt{m!}}
K(n,\kappa)^{\frac{m}{2}}\int_{\mathbb{R}^{n\times m}}^{'}
\frac{e^{i<x,z_{1}+\cdots+z_{m}> -\mu
t|z_{1}+\ldots+z_{m}|^{\alpha+\gamma}}} {(|z_{1}| \cdots
|z_{m}|)^{\frac{n-\kappa}{2}}}\overset{m}{\underset{l=1}{\prod}}W(dz_{l}),
\end{align}
where $\int^{'}\cdots$ denotes a $m$-fold Wiener integral with respect to
the complex Gaussian white noise $W(\cdot)$ on $\mathbb{R}^{n}$.
\end{Theorem}

\begin{bfseries}Remark.\end{bfseries}
To compare with Proposition 4 in the next Section 5, Theorem 1 has
the features that it involves both the Riesz and the Bessel
parameters, and that it also needs to rescale  the initial
condition. The intuitive meaning  behind the latter situation is
that, while the micro-scaling enforces to ``freeze down" both the
time $t$ and the space $x$, besides the overall renormalization we
also have to ``heat up" the initial data, in order to get a
non-degenerate (though singular) limiting field.
%%%%%%%%%%%%%%%%%%%%%%%%%%%%%%%% end of prop 4 %%%%%%%%%%%%%%%%%%%%%%%%%%%%%%%%%%%%%%%%%
%
%%%%%%%%%%%%%%%%%%%%%%%%%%% Begin of Theorem 2 %%%%%%%%%%%%%%%%%%%%%%%%%%%%%%%%%%%%%%%%%

Now, the micro-scaling of the system is

\begin{Theorem}\label{micro system theorem}
Let
$\mathbf{w}(t,x;\mathbf{w}_{0}(\cdot)):=(u(t,x;u_{0}(\cdot)),v(t,x;v_{0}(\cdot))),\
t> 0,\ x\in \mathbb{R}^{n}$, be the solution-vector of the initial
value problem (\ref{system}) and (\ref{initial}), satisfying the
Condition A, B, C and D. In the following, $\chi$ is a positive
parameter, $Q(t;d_{1},d_{2})$ is a matrix defined in (\ref{Q
matrix}), and  the two Gaussian noise fields $W_{j},\ j\in\{1,2\}$
are totally independent. Additionally, $m_{1},\ m_{2},\
\kappa_{1}\ and\ \kappa_{2}$ denote the parameters
contained in Condition C and D for $u_{0}$ and $v_{0}$.\\[0.2cm]
%%%%%%%%%%%% Theorem 2. Case 1.
%%%%%%%%%%%%%%%%%%%%%%%%%%%%%%%%%%%%%%%%%%%%%%%
(1) If\ \ $m_{2}\kappa_{2}>m_{1}\kappa_{1}$, then the
finite-dimensional distributions of the rescaled random field
\begin{align}\notag
[\varepsilon^{m_{1}\kappa_{1}\chi}L^{m_{1}}(\varepsilon^{-\chi})]^{-\frac{1}{2}}
\Big\{ \mathbf{w}(\varepsilon
t,\varepsilon^{\frac{1}{\alpha+\gamma}}x;\mathbf{w}_{0}(\varepsilon^{-\frac{1}{\alpha+\gamma}-\chi}\cdot))
-Q(\varepsilon
t;d_{1},d_{2})\left(\begin{array}{cccc}C^{(1)}_{0}\\C^{(2)}_{0}\end{array}\right)\Big\},\
\ t>0,\ x\in \mathbb{R}^{n},
\end{align}

\noindent converge weakly, as $\varepsilon \rightarrow 0$, to the
finite-dimensional distributions of the random field
\begin{align}\notag
\left(\begin{array}{cc}Y^{*}_{1}(t,x)\\Y^{*}_{2}(t,x)\end{array}\right):=\left(\begin{array}{cc}\widetilde{X}^{(1)}_{m_{1}}(t,x)\\
0 \end{array}\right),\ t>0,\ x\in \mathbb{R}^{n},
\end{align}
where
\begin{align}\label{limit field micro case 1}
\widetilde{X}^{(1)}_{m_{1}}(t,x):=\frac{C^{(1)}_{m_{1}}}{\sqrt{m_{1}!}}
K(n,\kappa_{1})^{\frac{m_{1}}{2}} \int_{\mathbb{R}^{n\times
m_{1}}}^{'} \frac{e^{i<x,z_{1}+\ldots+z_{m_{1}}> -\mu
t|z_{1}+\ldots+z_{m_{1}}|^{\alpha+\gamma}}} {(|z_{1}| \cdots
|z_{m_{1}}|)^{\frac{n-\kappa_{1}}{2}}}\prod_{l=1}^{m_{1}}{W_{1}(dz_{l})},
\end{align}
with $W_{1}(\cdot)$ is a complex Gaussian white noise on
$\mathbb{R}^{n}$ ( i.e., (\ref{limiting field(micro)}) with $m$,
$\kappa$ and $W$
replaced by $m_{1}$, $\kappa_{1}$ and $W_{1}$, respectively ).\\[0.2cm]
%%%%%%%%%%%%%%%%%%%%%%%% Theorem 2. Case 2. %%%%%%%%%%%%%%%%%%%%%%%%%%%%%%%%%%%%%%%%%%%%%%%
\noindent (2) If\ \ $m_{1}\kappa_{1}>m_{2}\kappa_{2}$, then the
finite-dimensional distributions of the rescaled random field
\begin{align}\notag
[\varepsilon^{m_{2}\kappa_{2}\chi}L^{m_{2}}(\varepsilon^{-\chi})]^{-\frac{1}{2}}
\Big\{ \mathbf{w}(\varepsilon
t,\varepsilon^{\frac{1}{\alpha+\gamma}}x;\mathbf{w}_{0}(\varepsilon^{-\frac{1}{\alpha+\gamma}-\chi}\cdot))
-Q(\varepsilon
t;d_{1},d_{2})\left(\begin{array}{cccc}C^{(1)}_{0}\\C^{(2)}_{0}\end{array}\right)\Big\},\
\ t>0,\ x\in \mathbb{R}^{n},
\end{align}
converge weakly, as $\varepsilon \rightarrow 0$, to the
finite-dimensional distributions of the random field
\begin{align}\notag
\left(\begin{array}{cc}Y^{**}_{1}(t,x)\\Y^{**}_{2}(t,x)\end{array}\right):=
\left(\begin{array}{cc}0
\\\widetilde{X}^{(2)}_{m_{2}}(t,x)\end{array}\right),\ t>0,\ x\in
\mathbb{R}^{n},
\end{align}
where
\begin{align}\label{limit field micro case 2}
\widetilde{X}^{(2)}_{m_{2}}(t,x):=\frac{C^{(2)}_{m_{2}}}{\sqrt{m_{2}!}}
K(n,\kappa_{2})^{\frac{m_{2}}{2}} \int_{\mathbb{R}^{n\times
m_{2}}}^{'} \frac{e^{i<x,z_{1}+\ldots+z_{m_{2}}> -\mu
t|z_{1}+\ldots+z_{m_{2}}|^{\alpha+\gamma}}} {(|z_{1}| \cdots
|z_{m_{2}}|)^{\frac{n-\kappa_{2}}{2}}}\prod_{l=1}^{m_{2}}W_{2}(dz_{l}),
\end{align}
and $W_{2}(\cdot)$ is a complex Gaussian white noise on
$\mathbb{R}^{n}$ ( i.e., (\ref{limiting field(micro)}) with $m$,
$\kappa$ and $W$
replaced by $m_{2}$, $\kappa_{2}$ and $W_{2}$, respectively ).\\[0.2cm]
%%%%%%%%%%%%%%%%%%%% Theorem 2. Case 3. %%%%%%%%%%%%%%%%%%%%%%%%%%%%%%%%%%%%%%%%%%%%%%%
\noindent(3) If\ \ $m_{1}=m_{2}:=m,\ \kappa_{1}=\kappa_{2}:=\kappa$,
then the finite-dimensional distributions of the rescaled random
field
\begin{align}\notag
[\varepsilon^{m\kappa\chi}L^{m}(\varepsilon^{-\chi})]^{-\frac{1}{2}}
\Big\{ \mathbf{w}(\varepsilon
t,\varepsilon^{\frac{1}{\alpha+\gamma}}x;\mathbf{w}_{0}(\varepsilon^{-\frac{1}{\alpha+\gamma}-\chi}\cdot))
-Q(\varepsilon
t;d_{1},d_{2})\left(\begin{array}{cccc}C^{(1)}_{0}\\C^{(2)}_{0}\end{array}\right)\Big\},\
\ t>0,\ x\in \mathbb{R}^{n},
\end{align}
converge weakly, as $\varepsilon \rightarrow 0$, to the
finite-dimensional distributions of the random field
\begin{align}\label{limit field micro case 3}
\left(\begin{array}{cc}Y^{***}_{1}(t,x)\\Y^{***}_{2}(t,x)\end{array}\right)
:=
\left(\begin{array}{cc}\widetilde{X}^{(1)}_{m}(t,x)\\\widetilde{X}^{(2)}_{m}(t,x)\end{array}\right)
,\ t>0,\ x\in \mathbb{R}^{n},
\end{align}
where $\widetilde{X}^{(1)}_{m}$ and  $\widetilde{X}^{(2)}_{m}$, are
defined in (\ref{limit field micro case 1}) and (\ref{limit field micro case 2})
with $m_{1}=m_{2}=m$ and $\kappa_{1}=\kappa_{2}=\kappa$.

\end{Theorem}

%%%%%%%%%%%%%%%%%%%%%%%%%%%%%%%%%%%%%%%%%%%%% End of Theorem 2 %%%%%%%%%%%%%%%%%%%%%%%%%%%%%%
%%%%%%%%%%%%%%%%%%%%%%%%%%%%%%%%%%%%%%%% Covariance structure for limit field %%%%%%%%%%%%%%%%%%%%%%%%
To understand the stochastic structure of the limiting fields, we
state, for instance, the following covariance result of
$(Y^{***}_{1}(t,x), Y^{***}_{2}(t,x))$.
\begin{Prop}\label{Cov matrix limit micro system}
For each fixed $t>0$, the limiting vector field
$\left(\begin{array}{cc}Y^{***}_{1}(t,x)&
Y^{***}_{2}(t,x)\end{array}\right)$ in the case (3) of Theorem {\ref{micro system theorem}} is
spatial-homogeneous and its covariance matrix has the following
spectral representation
\begin{align*}
E\left(\begin{array}{cc}Y^{***}_{1}(t,x)\\Y^{***}_{2}(t,x)\end{array}\right)
\left(\begin{array}{cc}Y^{***}_{1}(t^{'},x^{'}) &
Y^{***}_{2}(t^{'},x^{'})\end{array}\right)=
\int_{\mathbb{R}^{n}}e^{i<x-x^{'},\lambda>}S(\lambda;\alpha,\gamma)
d\lambda,
\end{align*}
where $ S(\lambda;\alpha,\gamma):=K(n,m\kappa) \frac{e^{-\mu (t+t^{'})
|\lambda|^{\alpha+\gamma}}}{(|\lambda|)^{n-m\kappa}} \left(
\begin{array}{cccc}
(C^{(1)}_{m})^{2} & 0 \\ 0 & (C^{(2)}_{m})^{2}
\end{array}
\right). $
\end{Prop}
%%%%%%%%%%%%%%%%%%%%%%%%%%%%%%%%%%%%%%%%%%%%%%%%%%%%%%%%%%%%%%%%%%%%%%%%%%%%%%%%%%%%%%%%%%%%
\begin{bfseries}Remark.\end{bfseries}
In view of the singularity of the spectral matrix near the origin,
we may conclude that, for limiting vector field in the case (3),
the long-range dependence (LRD) not only exists for each component
field but also exists between the two component fields; this is a
rather new phenomena for LRD, to our knowledge. Similar situation
happens for other cases, too.
%%%%%%%%%%%%%%%%%%%%%%%%%% 4. Limit theorems for solutions of P.D.E. systems %%%%%%%%%%%%%%%%%%%%%%%%%%%%%%%%%%%%
%%%%%%%%%%%%%%%%%%%%%%%%%%%%%%%%%%%%%%%%%%%%%%%%%%%%%%%%%%%%%%%%%%%%%%%%%%%
\section{Macro-scalings for the solution vector-field}
In this section, we present the macro-scaling limits for the
solution of the fractional kinetic systems  (\ref{system}) and (\ref{initial}),
in which only the Riesz parameter $\alpha$
plays its role in the scaling.\\
We again  begin with the following single-equation case, which is
adapted from \cite[Theorems 2.2 and 2.3]{AnhHomo}.
%%%%%%%%%%%%%%%%%%%%%%%%%% Proposition  (cite theorem). %%%%%%%%%%%%%%%%%%%%%%%%%%%%%%%
\begin{Prop} \label{cite theorem (Anh)}
Let
$s:=s(t,x;s_{0}(\cdot)),\ t> 0,\ x\in \mathbb{R}^{n}$, satisfies (\ref{single equation}), which
satisfies the above Condition B, C and D with $\kappa \in(0, n/m)$,
where $m$ denotes the Hermite rank
of the non-random function $h(\cdot)$ on $\mathbb{R}$,
which has the Hermite coefficients $C_{i}(h),\ i=0,1,\ldots$
(i.e., $h_{1}(\cdot)=h(\cdot),\ \zeta_{1}(x)=\zeta(x),\
f_{j}(\lambda)=f(\lambda)$ and $\kappa_{1}=\kappa$, etc.).
Then,\\
%%%%%%%%%%%%%%%%%%%%%
(1) The behaviour of the covariance function of the rescaled random
field $s^{\varepsilon}(t,x),\ t>0,\ x\in \mathbb{R}^{n}$,
\begin{align}\notag
s^{\varepsilon}(t,x):=[\varepsilon^{\frac{m\kappa}{\alpha}}L^{m}(\varepsilon^{-\frac{1}{\alpha}})]^{-\frac{1}{2}}
[s(\frac{t}{\varepsilon},\frac{x}{\varepsilon^{\frac{1}{\alpha}}};h(\zeta(\cdot)))-C_{0}(h)],\
t>0,\ x\in \mathbb{R}^{n}
\end{align}
is given by:
\begin{align}\label{citetheorem cov}
\underset{\varepsilon\rightarrow0}{\textup{lim}}\textup{Cov}(s^{\varepsilon}(t,x)s^{\varepsilon}(t^{'},x^{'}))
=(C_{m}(h))^{2}K(n,m\kappa)
\underset{\mathbb{R}^{n}}{\int}
\frac
{
e^{i<x-x^{'},\lambda>-\mu (t+t^{'})|\lambda|^{\alpha}}
}
{|\lambda|^{n-m\kappa}}
d\lambda.
\end{align}
\\
(2) Moreover, the finite dimensional distributions of the rescaled
field $s^{\varepsilon}(t,x)$ converge weakly, as $\varepsilon
\rightarrow 0$, to the finite-dimensional distributions of the
random field
\begin{align}\label{citetheorem}
s_{m}(t,x):=\frac{C_{m}(h)}{\sqrt{m!}}
K(n,\kappa)^{\frac{m}{2}}\int_{\mathbb{R}^{n\times m}}^{'}
\frac{e^{i<x,z_{1}+\cdots+z_{m}> -\mu
t|z_{1}+\ldots+z_{m}|^{\alpha}}} {(|z_{1}| \cdots
|z_{m}|)^{\frac{n-\kappa}{2}}}\overset{m}{\underset{l=1}{\prod}}W(dz_{l}),
\end{align}
$x\in \mathbb{R}^{n},\ t>0$, where
$\int^{'}\cdots$ denotes a $m$-fold Wiener integral with respect to the complex Gaussian white noise $W(\cdot)$ on $\mathbb{R}^{n}$.
\end{Prop}
%%%%%%%%%%%%%%%%%% End of Proposition 3. %%%%%%%%%%%%%%%%%%%%%%%%%%%%%%%%%%%%%%%%%
\begin{bfseries}Remark.\end{bfseries}
The above (\ref{citetheorem cov}) is expressed on the
``Fourier-domain ", which is more suitable for we will need; while
that (2.40) in \cite{AnhHomo} is in term of the variable domain.
%%%%%%%%%%%%%%%%%%%%%%%%%%%%% Theorem 2. %%%%%%%%%%%%%%%%%%%%%%%%%%%%%%%%%%

Then, the macro-scaling of the system is
\begin{Theorem}\label{macro system theorem}
Let
$\mathbf{w}(t,x;\mathbf{w}_{0}(\cdot)):=(u(t,x;u_{0}(\cdot)),v(t,x;v_{0}(\cdot))),\
t> 0,\ x\in \mathbb{R}^{n}$, be the solution-vector of the initial
value problem (\ref{system}) and (\ref{initial}), satisfying the
Condition A, B, C and D. In the following, $Q(t;d_{1},d_{2})$ is
the matrix defined in (\ref{Q matrix}), $p_{ij}$ is the entry in
(\ref{matrixB}), and the two Gaussian noise fields
$W_{j},\ j\in\{1,2\}$ are totally independent.
Additionally, $m_{1},\ m_{2},\ \kappa_{1}\ and\ \kappa_{2}$ denote the parameters
contained in Condition C and D for $u_{0}$ and $v_{0}$. \\[0.2cm]
%%%%%%%%%%%%%%%%%%%%%%%% Theorem 3 Case 1. %%%%%%%%%%%%%%%%%%%%%%%%%%%%%%%%%%%%%%%%%%%%%%%
(1) If\ \ $m_{2}\kappa_{2}>m_{1}\kappa_{1}$ and $d_{1}>d_{2}$, then
the finite-dimensional distributions of the rescaled random field
\begin{align}\notag
[\varepsilon^{\frac{m_{1}\kappa_{1}}{\alpha}}L^{m_{1}}(\varepsilon^{-\frac{1}{\alpha}})]^{-\frac{1}{2}}
e^{-d_{1}\frac{t}{\varepsilon}} \Big\{
\mathbf{w}(\frac{t}{\varepsilon},\frac{x}{\varepsilon^{\frac{1}{\alpha}}};\mathbf{w}_{0}(\cdot))
-Q(\frac{t}{\varepsilon};d_{1},d_{2})\left(\begin{array}{cccc}C^{(1)}_{0}\\C^{(2)}_{0}\end{array}\right)\Big\},\
\ t>0,\ x\in \mathbb{R}^{n},
\end{align}
converge weakly, as $\varepsilon \rightarrow 0$, to the
finite-dimensional distributions of the random field
\begin{align}\notag
\mathbf{T}_{m_{1}}^{(1)}(t,x)
:=\left(\begin{array}{cccc}p_{11}p_{22}X^{(1)}_{m_{1}}(t,x)\\p_{21}p_{22}X^{(1)}_{m_{1}}(t,x)\end{array}\right),\
\ t>0,\ x\in \mathbb{R}^{n},
\end{align}
where
\begin{align}\label{macro thm 1}
X^{(1)}_{m_{1}}(t,x):=\frac{C^{(1)}_{m_{1}}}{\sqrt{m_{1}!}}
K(n,\kappa_{1})^{\frac{m_{1}}{2}}\int_{\mathbb{R}^{n\times
m_{1}}}^{'} \frac{e^{i<x,z_{1}+\ldots+z_{m_{1}}> -\mu
t|z_{1}+\ldots+z_{m_{1}}|^{\alpha}}} {(|z_{1}| \cdots
|z_{m_{1}}|)^{\frac{n-\kappa_{1}}{2}}}\prod_{l=1}^{m_{1}}{W_{1}(dz_{l})},
\end{align}
with $W_{1}(\cdot)$ is a complex Gaussian white noise on
$\mathbb{R}^{n}$ ( i.e., (\ref{citetheorem}) with $m$, $\kappa$ and
$W$
replaced by $m_{1}$, $\kappa_{1}$ and $W_{1}$, respectively ).\\[0.2cm]
%%%%%%%%%%%%%%%%%%%%%%%%%% Theorem 3. Case 2. %%%%%%%%%%%%%%%%%%%%%%%%%%%%%%%%%%%%%%%%%%%%%%%
\noindent (2) If\ \ $m_{1}\kappa_{1}>m_{2}\kappa_{2}$ and
$d_{1}>d_{2}$, then the finite-dimensional distributions of the
rescaled random field
\begin{align}\notag
[\varepsilon^{\frac{m_{2}\kappa_{2}}{\alpha}}L^{m_{2}}(\varepsilon^{-\frac{1}{\alpha}})]^{-\frac{1}{2}}
e^{-d_{1}\frac{t}{\varepsilon}} \Big\{
\mathbf{w}(\frac{t}{\varepsilon},\frac{x}{\varepsilon^{\frac{1}{\alpha}}};\mathbf{w}_{0}(\cdot))
-Q(\frac{t}{\varepsilon};d_{1},d_{2})\left(\begin{array}{cccc}C^{(1)}_{0}\\C^{(2)}_{0}\end{array}\right)\Big\},\
\ t>0,\ x\in \mathbb{R}^{n},
\end{align}
converge weakly, as $\varepsilon \rightarrow 0$, to the
finite-dimensional distributions of the random field
\begin{align}\notag
\mathbf{T}_{m_{2}}^{(2)}(t,x)
:=\left(\begin{array}{cccc}-p_{11}p_{12}X^{(2)}_{m_{2}}(t,x)\\-p_{21}p_{12}X^{(2)}_{m_{2}}(t,x)\end{array}\right),\
\ t>0,\ x\in \mathbb{R}^{n},
\end{align}
where
\begin{align}\label{macro thm 2}
X^{(2)}_{m_{2}}(t,x):=\frac{C^{(2)}_{m_{2}}}{\sqrt{m_{2}!}}
K(n,\kappa_{2})^{\frac{m_{2}}{2}}\int_{\mathbb{R}^{n\times
m_{2}}}^{'} \frac{e^{i<x,z_{1}+\ldots+z_{m_{2}}> -\mu
t|z_{1}+\ldots+z_{m_{2}}|^{\alpha}}} {(|z_{1}| \cdots
|z_{m_{2}}|)^{\frac{n-\kappa_{2}}{2}}}\prod_{l=1}^{m_{2}}W_{2}(dz_{l}),
\end{align}
and $W_{2}(\cdot)$ is a complex Gaussian white noise on
$\mathbb{R}^{n}$ ( i.e., (\ref{citetheorem}) with $m$, $\kappa$ and
$W$
replaced by $m_{2}$, $\kappa_{2}$ and $W_{2}$, respectively ).\\[0.2cm]
%%%%%%%%%%%%%%%%%%%%%%%% Theorem 3 Case 3. %%%%%%%%%%%%%%%%%%%%%%%%%%%%%%%%%%%%%%%%%%%%%%%
(3) If\ \ $m_{1}=m_{2}:=m,\ \kappa_{1}=\kappa_{2}:=\kappa$,
and $d_{1}>d_{2}$, then the finite-dimensional distributions of the
rescaled random field
\begin{align}\notag
[\varepsilon^{\frac{m\kappa}{\alpha}}L^{m}(\varepsilon^{-\frac{1}{\alpha}})]^{-\frac{1}{2}}e^{-d_{1}\frac{t}{\varepsilon}}
\Big\{
\mathbf{w}(\frac{t}{\varepsilon},\frac{x}{\varepsilon^{\frac{1}{\alpha}}};\mathbf{w}_{0}(\cdot))
-Q(\frac{t}{\varepsilon};d_{1},d_{2})\left(\begin{array}{cccc}C^{(1)}_{0}\\C^{(2)}_{0}\end{array}\right)\Big\},\
\ t>0,\ x\in \mathbb{R}^{n},
\end{align}
converge weakly, as $\varepsilon \rightarrow 0$, to the
finite-dimensional distributions of the random field
\begin{align}
\mathbf{T}_{m}^{(3)}(t,x) := \mathbf{T}_{m}^{(1)}(t,x) +
\mathbf{T}_{m}^{(2)}(t,x),
\end{align}
where $\mathbf{T}_{m}^{(1)}(t,x)\  \textup{and}\
\mathbf{T}_{m}^{(2)}(t,x)$, are defined in the case (1)
and the case (2)
with $m_{1}=m_{2}=m$ and $\kappa_{1}=\kappa_{2}=\kappa$.\\
%
%%%%%%%%%%%%%%%%%%%%%%%% Theorem 3 case 4 %%%%%%%%%%%%%%%%%%%%%%%%%%%%%%%%%
\end{Theorem}
%%%%%%%%%%%%%%%%%%%%%%%%%%%%%%%%%%%%%
%%%%%%%%%%%%%%%%%%%%%%%%%%%%%%%%%%%%%

{\bf Remark:} In the above, we assume that $d_{1}>d_{2}$. In case
$d_{1}<d_{2}$, all the corresponding assertions hold, by
interchanging the roles of $m_1, \kappa_1$ and $m_2,\kappa_2$,
{\it etc.} As for $d_{1}=d_{2}$, it is reduced to the uncoupled
case and the result is induced from Proposition 4 directly.

\section{Time-fractional systems }
We extend the above results to the time-fractional derivative
$\frac{\partial^{\beta}}{\partial t^{\beta}},\ \beta\in(0,1),$ in
the system (\ref{system}), that is,
\begin{align}\label{extended system}
 \frac{\partial^{\beta}}{\partial t^{\beta}} \left(\begin{array}{cccc}u\\v\end{array}\right)=
 -\mu(I-\Delta)^{\frac{\gamma}{2}}
(-\Delta)^{\frac{\alpha}{2}}\left(\begin{array}{cccc}u\\v\end{array}\right)+
B\left(\begin{array}{cccc}u\\v\end{array}\right),\ \mu,\ \alpha,\
\gamma>0,
\end{align}
We recall that  the time-fractional derivative
$\frac{\partial^{\beta}}{\partial t^{\beta}}$ is  defined (see,
for example, the book of Djrbashian \cite{Djr})  by, for any
$\beta >0$,
\begin{align}\label{fraction def}
\frac{d^{\beta}f}{d
t^{\beta}}(t)=
\left\{
\begin{array}{lll}
f^{(m)}(t)  & \textup{if}\ \beta=m\in {\mathbb{N}}
\\
\frac{1}{\Gamma(m-\beta)}\int_{0}^{t}\frac{f^{(m)}(\tau)}{(t-\tau)^{\beta+1-m}}\
& \textup{if}\ \beta\in (m-1,m),
\end{array}
\right.\end{align} where $f^{(m)}(t)$ denotes the ordinary
derivative of order $m$ of a causal function $f(t)$ (i.e., $f$ is
vanishing for $t<0$).

The solution of (6.1) can be obtained, under Condition A, by
applying the Laplace and the Fourier transforms (see, for example,
\cite{Main,MainPara}), as follows.
\begin{align}\label{solution form of time fractional system}
{\bf{w}}(t,x;{\bf{w_{0}}}(\cdot))=
\int_{\mathbb{R}^{n}}P
\begin{pmatrix}G_{\beta}(t,x-y;d_{1})
& 0
\\0 & G_{\beta}(t,x-y;d_{2})
\end{pmatrix}
P^{-1}
\left(\begin{array}{cc}u_{0}(y)\\v_{0}(y)\end{array}\right)dy,
\end{align}
with the fractional Green function $G_{\beta}(t,x;d_{j})$ is defined by the
transformation
\begin{align*}
E_{\beta}(-\mu|\lambda|^{\alpha}(1+|\lambda|^{2})^{\frac{\gamma}{2}}t^{\beta}+d_{j}t^{\beta})
=\int_{\mathbb{R}^{n}}e^{i<x,\lambda>} G_{\beta}(t,x;d_{j})dx,\ \
j\in\{1,2\},
\end{align*}
where $E_{\beta}(\cdot)$ is the Mittag-Leffler function defined by
(see, for example, \cite{AnhHomo} or \cite[Chapter 1]{Djr})
\begin{align}
E_{\beta}(z)=\sum_{p=0}^{\infty}{\frac{z^{p}}{\Gamma (\beta p+1)}},\
z\in \mathbb{C},
\end{align}
and we shall use the following basic properties about the
Mittag-Leffler functions: they are entire functions on the complex
plane and their asymptotic behavior, when $\beta\in (0,2], \ \
\beta\neq 1,2$, has the inverse power law as follows:
\begin{align}\label{asymptoptic mittag}
|E_{\beta,\gamma}(z)|\ \sim\ O(\frac{1}{|z|}),\ \ |z|\
\rightarrow\ \infty\ \textup{with}\
|\textup{arg}(-z)|<\pi(1-\frac{\beta}{2}),\ \forall\ \gamma>0,
\end{align}
where arg: $\mathbb{C}\rightarrow(-\pi,\pi)$ and the notation
$f(z)\sim O(g(z))$ means that $f(z)/g(z)$ remains bounded as $z$
approaches the indicated limit point; see, for example, the
classic book by  Erd\'{e}lyi {\it et.al.} \cite{Erd} (pp. 206-212,
in particular p. 206 (7) and p. 210 (21)).

%which has the inverse power law (see, for example, Erd\'{e}lyi
%{\it
%et.al.} \cite{Erd} (pp.206-212, in particular p.206 (7) and p.210
%(21)))
%\begin{align}\label{asymptoptic mittag}
%|E_{\beta}(z)|\ \sim\ O(\frac{1}{|z|}),\ \ |z|\ \rightarrow\ \infty\
%\textup{with}\ |\textup{arg}(-z)|<\pi(1-\frac{\beta}{2}),\ \forall\
%\beta\in (0,1).
%\end{align}
We note that, from  \cite[(45)]{Wyss}, for $\beta\in(0,1)$, there
is an another representation for the fractional Green function
\begin{align}\label{Greenrelation}
G_{\beta}(t,x;d_{j})=t^{-\beta}\int^{\infty}_{0}f_{\beta}(t^{-\beta}s)G(s,x;d_{j})ds,\ t>0,\ x\in \mathbb{R}^{n},
\end{align}
where
\begin{align}\label{newheatkernel}
G(s,x;d_{j})=(\frac{1}{4\pi\mu s})^{\frac{n}{2}}e^{-\frac{|x|^{2}}{4\mu s}}e^{d_{j}s},\ s>0,\ x\in \mathbb{R}^{n},
\end{align}
while $f(p),p\geq 0,$ is a probability density which can be
represented by the H-function (see, for example, \cite[Section
3]{Wyss} and \cite[p. 284]{Schneider}) and its Laplace transform
is given by
\begin{align}\label{generater function}
\int_{0}^{\infty}e^{-qs}f_{\beta}(s)ds=E_{\beta}(-q),\ q\geq0.
\end{align}
Hence,
\begin{align}\notag
\int_{\mathbb{R}^{n}}G_{\beta}(t,x;d_{j})dx
\overset{(\ref{Greenrelation})}{=}
&\int_{\mathbb{R}^{n}}t^{-\beta}\int^{\infty}_{0}f_{\beta}(t^{-\beta}s)G(s,x;d_{j})ds\ dx\\\notag
=&t^{-\beta}\int_{0}^{\infty}f_{\beta}(t^{-\beta}s)
\int_{\mathbb{R}^{n}}G(s,x;d_{j})dx\ ds, (\textup{by Tonelli theorem})\\\label{integrable 1}
\overset{(\ref{newheatkernel})}{=}&t^{-\beta}\int_{0}^{\infty}f_{\beta}(t^{-\beta}s)e^{d_{j}s}\ ds
\overset{(\ref{generater function})}{=}E_{\beta}(d_{j}t^{\beta}),
\end{align}
where the convergence of the integral in (\ref{integrable 1}) is
guaranteed by the asymptotic behavior of the H-function (see, for
example, \cite[(3.7)]{Schneider}).

From the above discussion we know
$G_{\beta}(t,x;d_{j})\in L^{1}(\mathbb{R}^{n},dx)$ and
$\int_{\mathbb{R}^{n}}G_{\beta}(t,x;d_{j})dx=E_{\beta}(d_{j}t^{\beta})$
for any $t>0$ and $j\in\{1,2\}$.
Therefore, if the initial data $u_{0}(x)$ and $v_{0}(x)$ satisfy the
special form (\ref{initial data form}), then by the representation
(\ref{solution form of time fractional system}) we have
\begin{align}\label{expection matrix}
\mathbf{C}(t;B):=E{\bf{w}}(t,x;{\bf{w_{0}}}(\cdot))=
P\left(\begin{array}{cccc} E_{\beta}(d_{1}t^{\beta}) & 0
\\0
& E_{\beta}(d_{2}t^{\beta})
\end{array}\right)P^{-1}
\left(\begin{array}{cccc}C^{(1)}_{0}\\C^{(2)}_{0}\end{array}\right),
\end{align}
where $C^{(j)}_{0},\ j\in\{1,2\}$ are the Hermite coefficients
defined in (\ref{hermitecoeff}).

The following Theorems 4 and 5 are time-fractional versions those
in Theorems 3 and 2, respectively. However, it is now needed to
add the feature of the sub-diffusive property, which is the
reflection of time-fractional $\beta <1$ (see Section 1), into
consideration of the macro-scaling of homogenization. We need to
take an additional scaling on the matrix $B$ in the system
(\ref{extended system}) in order to compromise the effect of this
sub-diffusivity upon the interaction between $u$ and $v$. To
emphasize this situation, we denote the vector solution
$\mathbf{w}(t,x;\mathbf{w}_{0}(\cdot))$ by
$\mathbf{w}(t,x;\mathbf{w}_{0}(\cdot),B)$ in the following
formulation of macro-scaling of homogenization of a
spatial-temporal fractional kinetic system.

%%%%%%%%%%%%%%%%%%%% Section : Time fractional %%%%%%%%%%%%%%%%%%%%%%%%%%%%%%%%%%%%%%%%%%%%%
\begin{Theorem}
Let $\{{\bf{w}}(t,x;\mathbf{w}_{0}(\cdot),B),\ t>0,\ x\in
\mathbf{\mathbb{R}}^{n}\}$ be the solution-vector of the initial
value problem (\ref{extended system}) and (\ref{initial}),
satisfying Condition A, B, C and D. Moreover, the LRD parameter
$\kappa_{j}$ and the Hermite rank $m_{j}$ satisfy $m_{j}\kappa_{j} <
\textup{min}\{2\alpha,n\}$ and the Gaussian noise fields $W_{(j)}$
are totally independent for $j\in\{1,2\}$.
 \\[0.2cm]
%%%%%%%%%%%%%%%%%%%%%%%%%%%%%%%%%%%%%%%%%%%%%%%%% extend Theorem (macro). Case 1. %%%%%%%%%%%%%%%%%%%%%%%%%%%%%%%%%%%%%%%%%%%%%%%
(1) If $ m_{1}\kappa_{1}<m_{2}\kappa_{2}$, then the
finite-dimensional distributions of the rescaled random field
\begin{align}\notag
T_{\varepsilon}^{(1)}(t,x):=
(\varepsilon^{\frac{m_{1}\kappa_{1}}{\alpha}}
L^{m_{1}}(\varepsilon^{-\frac{\beta}{\alpha}}))^{-\frac{1}{2}}
\Big\{
{\bf{w}}(\varepsilon^{-1}t,\varepsilon^{-\frac{\beta}{\alpha}}x;\mathbf{w}_{0}(\cdot),\varepsilon^{\beta}B)
-\mathbf{C}(\varepsilon^{-1}t;\varepsilon^{\beta}B) \Big\},
\end{align}
$\ t>0,\ x\in \mathbb{R}^{n},$ converge weakly, as $\varepsilon
\rightarrow 0$, to the finite-dimensional distributions of the
random field
\begin{align}\notag
T^{(1)}(t,x)= \left(\begin{array}{cccc}
p_{11}p_{22}T^{(1)}(t,x;d_{1}) -p_{12}p_{21}T^{(1)}(t,x;d_{2})
\\
p_{21}p_{22}T^{(1)}(t,x;d_{1}) -p_{21}p_{22}T^{(1)}(t,x;d_{2})
\end{array}\right),\ t>0,\ x \in \mathbb{R}^{n},
\end{align}
where for $j\in\{1,2\}$
\begin{align}\label{fractional time limiting field 1}
&T^{(1)}(t,x;d_{j})
\\\notag
:=& \frac{C^{(1)}_{m_{1}}}{\sqrt{m_{1}!}}
K(n,\kappa_{1})^{\frac{m_{1}}{2}} \underset{\mathbb{R}^{n\times
m_{1}}}{\int^{'}} e^{i<x,\lambda_{1}+\cdots+\lambda_{m_{1}}>}
\frac{E_{\beta}(-\mu|\lambda_{1}+\cdots+\lambda_{m_{1}}|^{\alpha}t^{\beta}+d_{j}t^{\beta})}
{(|\lambda_{1}| \cdots |\lambda_{m_{1}}|)^{\frac{n-\kappa_{1}}{2}}}
\overset{m_{1}}{\underset{l=1}{\prod}}W_{1}(d\lambda_{l}).
\end{align}
%
%
%
%%%%%%%%%%%%%%%%%%%%%%%%%%%%%%%%%%%%%%%%%%%% Extend Theorem macro . Case 2. %%%%%%%%%%%%%%%%%%%%%%%%%%%%%%%%%
(2) If $ m_{2}\kappa_{2}<m_{1}\kappa_{1}$, then the
finite-dimensional distributions of the rescaled random field
\begin{align}\notag
T_{\varepsilon}^{(2)}(t,x)
(\varepsilon^{\frac{m_{2}\kappa_{2}}{\alpha}}
L^{m_{2}}(\varepsilon^{-\frac{\beta}{\alpha}}))^{-\frac{1}{2}}
\Big\{
{\bf{w}}(\varepsilon^{-1}t,\varepsilon^{-\frac{\beta}{\alpha}}x;\mathbf{w}_{0}(\cdot),\varepsilon^{\beta}B)
-\mathbf{C}(\varepsilon^{-1}t;\varepsilon^{\beta}B) \Big\},
\end{align}
$\ t>0,\ x\in \mathbb{R}^{n},$ converge weakly, as $\varepsilon
\rightarrow 0$, to the finite-dimensional distributions of the
random field
\begin{align}\notag
T^{(2)}(t,x)= \left(\begin{array}{cccc}
-p_{11}p_{12}T^{(2)}(t,x;d_{1})+p_{11}p_{12}T^{(2)}(t,x;d_{2})
\\
-p_{12}p_{21}T^{(2)}(t,x;d_{1})+p_{11}p_{22}T^{(2)}(t,x;d_{2})
\end{array}\right),\ t>0,\ x \in \mathbb{R}^{n},
\end{align}
where for $j\in\{1,2\}$
\begin{align}\label{fractional time limiting field 2}
&T^{(2)}(t,x;d_{j})
\\\notag
:= &\frac{C^{(2)}_{m_{2}}}{\sqrt{m_{2}!}}
K(n,\kappa_{2})^{\frac{m_{2}}{2}} \underset{\mathbb{R}^{n\times
m_{2}}}{\int^{'}} e^{i<x,\lambda_{1}+\cdots+\lambda_{m_{2}}>}
\frac{E_{\beta}(-\mu|\lambda_{1}+\cdots+\lambda_{m_{2}}|^{\alpha}t^{\beta}+d_{j}t^{\beta})}
{(|\lambda_{1}| \cdots |\lambda_{m_{2}}|)^{\frac{n-\kappa_{2}}{2}}}
\overset{m_{2}}{\underset{l=1}{\prod}}W_{2}(d\lambda_{l}).
\end{align}
%
%
%%%%%%%%%%%%%%%%%% Extended Theorem macro . Case 3. %%%%%%%%%%%%%%%%%%%%%%%%%%%%%%%%%%%%
(3) If $ m_{1}=m_{2}=m$ and $\kappa_{1}=\kappa_{2}=\kappa$, then the
finite-dimensional distributions of the rescaled random field
\begin{align}\notag
T_{\varepsilon}^{(3)}(t,x):= (\varepsilon^{\frac{m\kappa}{\alpha}}
L^{m}(\varepsilon^{-\frac{\beta}{\alpha}}))^{-\frac{1}{2}} \Big\{
{\bf{w}}(\varepsilon^{-1}t,\varepsilon^{-\frac{\beta}{\alpha}}x;\mathbf{w}_{0}(\cdot),\varepsilon^{\beta}B)
-\mathbf{C}(\varepsilon^{-1}t;\varepsilon^{\beta}B) \Big\},
\end{align}
$\ t>0,\ x\in \mathbb{R}^{n},$ converge weakly, as $\varepsilon
\rightarrow 0$, to the finite-dimensional distributions of the
random field
\begin{align}\notag
T^{(3)}(t,x)=\widetilde{T}^{(1)}(t,x)+\widetilde{T}^{(2)}(t,x),\
t>0,\ x \in \mathbb{R}^{n},
\end{align}
where the random field $\widetilde{T}^{(j)}(t,x)$ is the same as
the limiting random field $T^{(j)}(t,x)$ in (\ref{fractional time
limiting field 1}) and (\ref{fractional time limiting field 2}) by
replacing $m_{j}\rightarrow m$ and $\kappa_{j}\rightarrow \kappa$
for $j\in \{1,2\}$.
\end{Theorem}
%%%%%%%%%%%%%%%%%%%% End of  extended Theorem (macro) . %%%%%%%%%%%%%%%%%%%%%%%%%%%%%%%%%%%%%%%%%%%%%%%
{\bf Remark.} The restriction
$m_{j}\kappa_{j}<\textup{min}\{2\alpha,n\}$ in the above Theorem
4, together with the power law decay of Mittag-Leffler functions,
guarantee that the random fields $T^{(j)}(t,x),\ j\in\{1,2,3\}$,
are indeed defined as $L^{2}(\Omega, \mathcal{F}, \mathcal{P})$
stochastic integrals.
%

%%%%%%%%%%%%%%%%%%%%%%%%%%%%%%%%%%%%%%%%%%%%%%%%%%%%%%%%%%%%%%%%%%%%%%%%%%%%%%%%%%%%%%%%%%%%
%%%%%%%%%%%%%%%%%%%%%%%%%%% Begin  the discussion of Micro %%%%%%%%%%%%%%%%%%%%%%%%%%%%%%%%%%%%%%%%
As for the micro-scaling, the sub-diffusivity has no influence,
and the same micro-scaling procedure as Theorem 2 applies.
\begin{Theorem}
Let $\{{\bf{w}}(t,x;\mathbf{w}_{0}(\cdot),B),\ t>0,\ x\in
\mathbf{\mathbb{R}}^{n}\}$ be the solution-vector of the initial
value problem (\ref{extended system}) and (\ref{initial}),
satisfying Condition A, B, C and D. Moreover, the LRD parameter
$\kappa_{j}$ and the Hermite rank $m_{j}$ satisfy
$m_{j}\kappa_{j}<\textup{min}\{2(\alpha+\gamma),n\}$ and the
Gaussian noise fields $W_{(j)}$ are totally independent for
$j\in\{1,2\}$.
 \\[0.2cm]
%%%%%%%%%%%%%%%%%%%%%%%%% extend Theorem (micro). Case 1. %%%%%%%%%%%%%%%%%%%%%%%%%%%%%%%%%%%%%%%%%%%%%%%
(1) If $ m_{1}\kappa_{1}<m_{2}\kappa_{2}$, then the
finite-dimensional distributions of the rescaled random field
\begin{align}\notag
\mathbf{M}^{(1)}_{\varepsilon}:=[\varepsilon^{m_{1}\kappa_{1}\chi}L^{m_{1}}(\varepsilon^{-\chi})]^{-\frac{1}{2}}
\Big\{ \mathbf{w}(\varepsilon
t,\varepsilon^{\frac{\beta}{\alpha+\gamma}}x;\mathbf{w}_{0}(\varepsilon^{-\frac{\beta}{\alpha+\gamma}-\chi}\cdot))
-\left(\begin{array}{cc}C^{(1)}_{0}\\
C^{(2)}_{0}\end{array}\right)\Big\},\ t>0,\ x\in \mathbb{R}^{n},
\end{align}
converge weakly, as $\varepsilon \rightarrow 0$, to the
finite-dimensional distributions of the random field
\begin{align}\notag
\mathbf{M}^{(1)}(t,x)= \left(\begin{array}{cccc} M^{(1)}(t,x)
\\
0
\end{array}\right),\ t>0,\ x \in \mathbb{R}^{n},
\end{align}
where
\begin{align*}
M^{(1)}:= \frac{C^{(1)}_{m_{1}}}{\sqrt{m_{1}!}}
K(n,\kappa_{1})^{\frac{m_{1}}{2}} \underset{\mathbb{R}^{n\times
m_{1}}}{\int^{'}} e^{i<x,\lambda_{1}+\cdots+\lambda_{m_{1}}>}
\frac{E_{\beta}(-\mu|\lambda_{1}+\cdots+\lambda_{m_{1}}|^{\alpha+\gamma}t^{\beta})}
{(|\lambda_{1}| \cdots |\lambda_{m_{1}}|)^{\frac{n-\kappa_{1}}{2}}}
\overset{m_{1}}{\underset{l=1}{\prod}}W_{1}(d\lambda_{l}).
\end{align*}
%
%
%
%%%%%%%%%%%%%%%%%%%%%%%%%%%%%%%%%%%%%%%%%%%% Extend Theorem micro . Case 2. %%%%%%%%%%%%%%%%%%%%%%%%%%%%%%%%%
(2) If $ m_{2}\kappa_{2}<m_{1}\kappa_{1}$, then the
finite-dimensional distributions of the rescaled random field
\begin{align}\notag
\mathbf{M}^{(2)}_{\varepsilon}:=[\varepsilon^{m_{2}\kappa_{2}\chi}L^{m_{2}}(\varepsilon^{-\chi})]^{-\frac{1}{2}}
\Big\{ \mathbf{w}(\varepsilon
t,\varepsilon^{\frac{\beta}{\alpha+\gamma}}x;\mathbf{w}_{0}(\varepsilon^{-\frac{\beta}{\alpha+\gamma}-\chi}\cdot))
-\left(\begin{array}{cc}C^{(1)}_{0}\\
C^{(2)}_{0}\end{array}\right)\Big\},\ t>0,\ x\in \mathbb{R}^{n},
\end{align}
 converge weakly, as $\varepsilon
\rightarrow 0$, to the finite-dimensional distributions of the
random field
\begin{align}\notag
\mathbf{M}^{(2)}(t,x)= \left(\begin{array}{cccc} 0
\\
M^{(2)}(t,x)
\end{array}\right),\ t>0,\ x \in \mathbb{R}^{n},
\end{align}
where
\begin{align*}
M^{(2)}:= \frac{C^{(2)}_{m_{2}}}{\sqrt{m_{2}!}}
K(n,\kappa_{2})^{\frac{m_{2}}{2}} \underset{\mathbb{R}^{n\times
m_{2}}}{\int^{'}} e^{i<x,\lambda_{1}+\cdots+\lambda_{m_{2}}>}
\frac{E_{\beta}(-\mu|\lambda_{1}+\cdots+\lambda_{m_{2}}|^{\alpha+\gamma}t^{\beta})}
{(|\lambda_{1}| \cdots |\lambda_{m_{2}}|)^{\frac{n-\kappa_{2}}{2}}}
\overset{m_{2}}{\underset{l=1}{\prod}}W_{2}(d\lambda_{l}).
\end{align*}
%
%
%%%%%%%%%%%%%%%%%% Extended Theorem micro . Case 3. %%%%%%%%%%%%%%%%%%%%%%%%%%%%%%%%%%%%
(3) If $ m_{1}=m_{2}=m$ and $\kappa_{1}=\kappa_{2}=\kappa$, then the
finite-dimensional distributions of the rescaled random field
\begin{align}\notag
\mathbf{M}^{(3)}_{\varepsilon}:=[\varepsilon^{m\kappa\chi}L^{m}(\varepsilon^{-\chi})]^{-\frac{1}{2}}
\Big\{ \mathbf{w}(\varepsilon
t,\varepsilon^{\frac{\beta}{\alpha+\gamma}}x;\mathbf{w}_{0}(\varepsilon^{-\frac{\beta}{\alpha+\gamma}-\chi}\cdot))
-\left(\begin{array}{cc}C^{(1)}_{0}\\
C^{(2)}_{0}\end{array}\right)\Big\},\ t>0,\ x\in \mathbb{R}^{n},
\end{align}
converge weakly, as $\varepsilon \rightarrow 0$, to the
finite-dimensional distributions of the random field
\begin{align}\notag
\mathbf{M}^{(3)}(t,x)=\widetilde{\mathbf{M}}^{(1)}(t,x)+\widetilde{\mathbf{M}}^{(2)}(t,x),\
t>0,\ x \in \mathbb{R}^{n},
\end{align}
where the random field $\widetilde{\mathbf{M}}^{(j)}(t,x)$ is the
same as the limiting random field $\mathbf{M}^{(j)}(t,x)$ of Case
(1) and (2) by replacing $m_{j}\rightarrow m$ and
$\kappa_{j}\rightarrow \kappa$ for $j\in \{1,2\}$.
\end{Theorem}
{\bf The concluding remark:} The time-fractional index $\beta < 1$
indicates the sub-diffusivity, and it changes to be the
super-diffusivity if we consider $\beta >1$ (see Section 1). In
\cite{LSh}, the time-fractional reaction-wave type system with
random initial data are studied, in which the first-order
time-derivatives of the initial data play the crucial role. To
consider spatial-temporal fractional kinetic systems which is
super-diffusive in time and Riesz-Bessel in space will be a task
of tremendous analysis. Finally, we mention that, for  the
classical, {\it i.e.} non-fractional, heat-type system with random
initial condition, the solution vector-field and the scaling limit
are expressed in terms of heat kernels; this more explicit and
simpler case is treated in \cite{LSh08}.
%%%%%%%%%%%%%%%%%%%%%%%%%%%%%%%%%%%%  End of extend theorem (micro)   %%%%%%%%%%%%%%%%%%%%%%%%%%%
%%%%%%%%%%%%%%%%%%%%%%%%%%%%%%%%%%%%%%%%%%%%%%%%%%%%%%%%%%%%%%%%%%%%%%%%%%%%%%%%%%%%%%%%%%%%%%%%%%%%%

\vskip 20 pt
%%%%%%%%%%%%%%%%%%%%%%%%%%%  Proof %%%%%%%%%%%%%%%%%%%%%%%%%%%%%%%%%%%%%%

\section{Proofs}

In the following proofs, $\Rightarrow$\ \  denotes the convergence
of random variables (or random families) in distributional sense,
and $\overset{d}{=}$\ \ denotes the equality of random variables (or random families) in distributional sense.
Moreover, we also denote $f(t,x;\varepsilon)\asymp g(t,x;\varepsilon)$
if there exists a constant $c:=c(t,x)>0$ such that $cg(t,x;\varepsilon)<f(t,x;\varepsilon)<c^{-1}g(t,x;\varepsilon)$
when $\varepsilon\rightarrow 0$.\\
%
%%%%%%%%%%%%%%%% Proof of Proposition 1. %%%%%%%%
             %%%%%%%%%%%%%%%%%%%%%%%%%%%%%%%%%%%%
{\bf Proof of Proposition 1.}\\
For (\ref{spectral solution}), we use the solution form (\ref{solution}) and Karhunen's
representation (\ref{Spectral represent}) to get
\begin{align}\label{proofpro2matrix}
{\bf{w}}(t,x;\mathbf{w}_{0}(\cdot))&=Q(t,d_{1},d_{2})
\int_{\mathbb{R}^{n}}G(t,y;\alpha,\gamma) \left(\begin{array}{cc}
\int_{\mathbb{R}^{n}}e^{i<\lambda,x-y>}Z_{F_{1}}(d\lambda)\\
\int_{\mathbb{R}^{n}}e^{i<\lambda,x-y>}Z_{F_{2}}(d\lambda)
\end{array}\right)dy\\\notag
&=Q(t,d_{1},d_{2})
\int_{\mathbb{R}^{n}}\int_{\mathbb{R}^{n}}e^{i<\lambda,x-y>}G(t,y;\alpha,\gamma)dy \left(\begin{array}{cc}
Z_{F_{1}}(d\lambda)\\
Z_{F_{2}}(d\lambda)
\end{array}\right)\\\notag
&\overset{(\ref{fractional
Green function})}{=}
Q(t,d_{1},d_{2})
\int_{\mathbb{R}^{n}}e^{i<\lambda,x>}e^{-\mu
t|\lambda|^{\alpha}(1+|\lambda|^{2})^{\frac{\gamma}{2}}}\left(\begin{array}{cc}
Z_{F_{1}}(d\lambda)\\
Z_{F_{2}}(d\lambda)
\end{array}\right)
.
\end{align}
For (\ref{covariance structure}), it is a consequence by using the
independence assumption between the initial data, and we proceed it
as follows,
\begin{align*}
&E{\bf{w}}(t,x;\mathbf{w}_{0}(\cdot))\overline{{\bf{w}}(t^{'},x^{'};\mathbf{w}_{0}(\cdot))}
\\
=& Q(t;d_{1},d_{2}) \int_{\mathbb{R}^{n}}
e^{i(<\lambda,x>-<\lambda^{'},x^{'}>)} e^{-\mu
t|\lambda|^{\alpha}(1+|\lambda|^{2})^{\frac{\gamma}{2}}-\mu
t^{'}|\lambda^{'}|^{\alpha}(1+|\lambda^{'}|^{2})^{\frac{\gamma}{2}}}\times
\\
&E\left(\begin{array}{cc}
Z_{F_{1}}(d\lambda)\\Z_{F_{2}}(d\lambda)\end{array}\right)
\left(\begin{array}{cc}
Z_{F_{1}}(d\lambda^{'})\\Z_{F_{2}}(d\lambda^{'})\end{array}\right)^{T}Q(t^{'};d_{1},d_{2})^{T}
\\
=& Q(t;d_{1},d_{2}) \int_{\mathbb{R}^{n}} e^{i(<\lambda,x-x^{'}>)}
e^{-\mu
(t+t^{'})|\lambda|^{\alpha}(1+|\lambda|^{2})^{\frac{\gamma}{2}}}
\left(\begin{array}{cccc} F_{1}(d\lambda) & 0\\ 0 & F_{2}(d\lambda)
\end{array}\right)
Q(t^{'};d_{1},d_{2})^{T}.\hspace{0.6cm}\qed
\end{align*}
%%%%%%%%%%% Slutsky's argument ... %%%%%%%%%%%%%%%%%%%%%%%%%%%%%

Before going to prove our main results we
recall the following two arguments, which are
powerful to help us to reduce and simplify our problems.\\
\textbf{(Slutsky argument)} (see, for example, the book of
Leonenko \cite[p. 6.]{L})  Let  $\{\xi_\varepsilon\}$ and
$\{\eta_{\varepsilon}\}$ be families of random variables such that
$\{\xi_{\varepsilon}\}\ \Rightarrow\ \xi$ and $\
\{\eta_{\varepsilon}\} \overset{P}{\rightarrow}\ c$, where $c\in
R$. Then $(i)\ \xi_{\varepsilon}+\eta_{\varepsilon}\ \Rightarrow
\xi\ +\ c$, $(ii)\ \xi_{\varepsilon}\eta_{\varepsilon}\
\Rightarrow\ c\xi$, and
$(iii)\ \xi_{\varepsilon}/\eta_{\varepsilon}\ \Rightarrow\ \xi/c.$\\[0.2cm]
%%%
\textbf{(Cramer-Wold argument)} (see, for example, again  \cite[p.
6.]{L})  A family of $k$-dimensional r.v's
$\xi_{\varepsilon}:=(\xi_{\varepsilon}(x_{1}),...,\xi_{\varepsilon}(x_{k}))^{T}$
converge in distribution to a $k$-dimensional r.v.
$\xi:=(\xi(x_{1}),...,\xi(x_{k}))^{T}$ if and only if, for every
$c:=(c_{1},...c_{k})^{T}\in \mathbb{R}^{k},$
\[<c,\xi_{\varepsilon}>=\sum^{k}_{i=1}c_{i}\xi_{\varepsilon}(x_{i})\Rightarrow \sum^{k}_{i=1}c_{i}\xi(x_{i})=<c,\xi>,\ \ \textup{as}\ \varepsilon \rightarrow0.\]

%%%%%%%%%%%%%%%%%%%%%%%%%%%%%%%%%%%%%%%%% Extension of Cramer Wold argument %%%%%%%%%%%%%%%%%%%%%%%%%%%%%%%%%%%%%%%%%%
The following lemma, although it is a simple extension of
Cramer-Wold argument, is of interest in itself, and will play an
important role in the proof
of our Theorem 2 and 3.\\[0.2cm]
\textbf{Lemma 1.}\ Let $X_{\varepsilon}:=\{[X_{\varepsilon}^{(1)},\
X_{\varepsilon}^{(2)}]^{T}(t,x),\ x\in \mathbb{R}^{n},\ t>0\}$ be a
$\mathbb{R}^{2}$-valued random field which is generated by
$X_{\varepsilon}(t,x)=Q_{\varepsilon}(t)[U_{\varepsilon},\
V_{\varepsilon}]^{T}(t,x)$, where $U_{\varepsilon}(t,x)$ and
$V_{\varepsilon}(t,x)$ are independent random fields on
$\mathbb{R}^{n}\times \mathbb{R}^{+}$ and $Q_{\varepsilon}(t)$ is a
non-random $2\times 2$ matrix. If there exist two random fields
$U_{0}$ and $V_{0}$ such that $U_{\varepsilon}(t,x)\Longrightarrow
U_{0}(t,x)$ and $V_{\varepsilon}(t,x)\Longrightarrow V_{0}(t,x)$,
respectively, and $Q_{\varepsilon}(t)$ converges to $Q_(t)$ in the
usual sense when $\varepsilon\rightarrow 0$, then the finite
dimensional distributions of $X_{\varepsilon}(t,x),\ t>0,\ x\in
\mathbb{R}^{n}$, converge to the finite dimensional distributions of
$X:=\{Q(t)[U_{0},\ V_{0}]^{T}(t,x),\ t>0,\ x\in \mathbb{R}^{n}\}.$\\[0.2cm]
%%%%%%%%%%%%%%%%%%%%%%%%%%%%%%%%%%%%% {\bf Proof of Lemma}%%%%%%%%%%%%%%%%%%%%
{\bf Proof of Lemma 1}\\
By Cramer-Wold argument with $k=2$ there, it suffices to prove: For any given
$c_{1}$, $c_{2}\in \mathbb{R}$ and $x$, $t$ fixed, we have
\begin{align}\notag
[c_{1},\ \ c_{2}][Q_{\varepsilon}(t)][U_{\varepsilon}(x,t),\ \
V_{\varepsilon}(x,t)]^{T}\Rightarrow [c_{1},\ \
c_{2}][Q_(t)][U_{0}(x,t),\ \  V_{0}(x,t)]^{T},
\end{align}
which is equivalent to\\
$
\left(c_{1}Q_{11,\varepsilon}(t)+c_{2}Q_{21,\varepsilon}(t)\right)U_{\varepsilon}(x,t)+
\left(c_{1}Q_{12,\varepsilon}(t)+c_{2}Q_{22,\varepsilon}(t)\right)V_{\varepsilon}(x,t)
\\
\Rightarrow \left(c_{1}Q_{11}(t)+c_{2}Q_{21}(t)\right)U_{0}+
\left(c_{1}Q_{12}(t)+c_{2}Q_{22}(t)\right)V_{0},$\\
where the $i,j$ indicate the $(i,j)$ entry of the matrix. While the
above display can be checked by using the characteristic functions,
since  $U_{\varepsilon}$ and $V_{\varepsilon}$
 are assumed to be independent(whence
so are $U_0,V_0$).\qed\\
%%%%%%%
%%%%%%
%%%%%%%%%%%%%%%%%%%%%%%% proof of (micro) single Theorem . %%%%%%%%%%%%%%%%%%%%%%%%%
                                %%%%%%%%%%%%%%%%%%%%%%%%%%%%%%%%%%%%%%%%%%%%%%%%%%%%%%%%%%%%
{\bf Proof of Theorem 1.}\\
%%%%%%%%%%%%%%%%%%%%%%%% Proposition 4 (1) %%%%%%%%%%%%%%%%%%%%%%%%%%
{\bf (1)} Firstly, for simplification, we set
$N_{1}(\varepsilon):=\varepsilon^{\chi m\kappa}L^{m}(\varepsilon^{-\chi})$.
By the Hermite expansion, we can rewrite
$s(\varepsilon t, \varepsilon^{\frac{1}{\alpha+\gamma}}x
;h(\zeta(\varepsilon^{-\frac{1}{\alpha+\gamma}-\chi}\cdot)))-C_{0}(h)$
as
\begin{align}\label{micro 1 decomposition}
s(\varepsilon t, \varepsilon^{\frac{1}{\alpha+\gamma}}x
;h(\zeta(\varepsilon^{-\frac{1}{\alpha+\gamma}-\chi}\cdot)))-C_{0}(h)
= \sum_{\rho=m}^{\infty}{s(\varepsilon
t,\varepsilon^{\frac{1}{\alpha+\gamma}}x;\frac{C_{\rho}(h)}{\sqrt{\rho!}}
H_{\rho}(\zeta(\varepsilon^{-\frac{1}{\alpha+\gamma}-\chi}\cdot)))},
\end{align}
where the summation is in $L^2(\Omega)$ sense.  Hence, in
accordion to the definition (\ref{rescaled micro field}) about the
random field $s^{\varepsilon}(t,x)$, it can be rewritten as
\begin{equation}\label{expansion single micro}
s^{\varepsilon}(t,x)=\overset{\infty}{\underset{\rho=m}{\sum}}I^{\varepsilon}_{\rho}(t,x),
\end{equation}
with
$I^{\varepsilon}_{\rho}(t,x):=(N_{1}(\varepsilon))^{-\frac{1}{2}}s(\varepsilon
t,\varepsilon^{\frac{1}{\alpha+\gamma}}x;\frac{C_{\rho}(h)}{\sqrt{\rho!}}
H_{\rho}(\zeta(\varepsilon^{-\frac{1}{\alpha+\gamma}-\chi}\cdot)))
.$\\
By (\ref{solution}) with $d_{1}=d_{2}=0$,
we have
\begin{align}\label{proof micro single add1}
&(\frac{C_{\rho}(h)}{\sqrt{\rho!}})^{-1}(N_{1}(\varepsilon))^{\frac{1}{2}}I^{\varepsilon}_{\rho}(t,x) \\\notag
&=
\underset{\mathbb{R}^{n}}{\int}G(\varepsilon
t,y;\alpha,\gamma)H_{\rho}(\zeta(\varepsilon^{-\frac{1}{\alpha+\gamma}-\chi}(\varepsilon^{\frac{1}{\alpha+\gamma}}x-y)))dy
\\\notag
&= \underset{\mathbb{R}^{n}}{\int}G(\varepsilon t,y;\alpha,\gamma)
\underset{\mathbb{R}^{n\times \rho}}{\int^{'}}
e^{i<\varepsilon^{-\frac{1}{\alpha+\gamma}-\chi}(\varepsilon^{\frac{1}{\alpha+\gamma}}x-y),\lambda_{1}+\cdots+\lambda_{\rho}>}
\overset{\rho}{\underset{{\sigma=1}}{\prod}}{\sqrt{f(\lambda_{\sigma})}}
W(d\lambda_{\sigma})dy
\\\notag
&=\underset{\mathbb{R}^{n\times \rho}}{\int^{'}}
e^{i<\varepsilon^{-\chi}x,\lambda_{1}+\cdots+\lambda_{\rho}>}
\Big\{
\underset{\mathbb{R}^{n}}{\int}G(\varepsilon t,y;\alpha,\gamma)
e^{-i<\varepsilon^{-\frac{1}{\alpha+\gamma}-\chi}y,\lambda_{1}+\cdots+\lambda_{\rho}>}
\Big\}
\overset{\rho}{\underset{{\sigma=1}}{\prod}}{\sqrt{f(\lambda_{\sigma})}}
W(d\lambda_{\sigma})dy.
\end{align}
For the bracket above, by substituting
$t\rightarrow \varepsilon t$
and
$\lambda\rightarrow \varepsilon^{-\frac{1}{\alpha+\gamma}-\chi}(\lambda_{1}+\cdots+\lambda_{\rho})$
into (\ref{fractional Green function}), we have
\begin{equation*}
\underset{\mathbb{R}^{n}}{\int}G(\varepsilon t,y;\alpha,\gamma)
e^{-i<\varepsilon^{-\frac{1}{\alpha+\gamma}-\chi}y,\lambda_{1}+\cdots+\lambda_{\rho}>}
=
e^{-\mu\varepsilon
t\varepsilon^{-\frac{\alpha}{\alpha+\gamma}-\alpha\chi}|\lambda_{1}+\cdots+\lambda_{\rho}|^{\alpha}
(1+\varepsilon^{-\frac{2}{\alpha+\gamma}-2\chi}|\lambda_{1}+\cdots+\lambda_{\rho}|^{2})^{\frac{\gamma}{2}}}
,
\end{equation*}
so (\ref{proof micro single add1}) is equal to
\begin{align}\notag
&\underset{{\mathbb{R}^{n\times \rho}}}{\int^{'}}
e^{i<\varepsilon^{-\chi}x,\lambda_{1}+\cdots+\lambda_{\rho}>}
e^{-\mu\varepsilon
t\varepsilon^{-\frac{\alpha}{\alpha+\gamma}-\alpha\chi}|\lambda_{1}+\cdots+\lambda_{\rho}|^{\alpha}
(1+\varepsilon^{-\frac{2}{\alpha+\gamma}-2\chi}|\lambda_{1}+\cdots+\lambda_{\rho}|^{2})^{\frac{\gamma}{2}}
}
\prod_{\sigma=1}^{\rho}{\sqrt{f(\lambda_{\sigma})}}
W(d\lambda_{\sigma})\\\label{proof micro single 1}
&\overset{d}{=} \varepsilon^{\frac{\chi\rho n}{2}}
\underset{{\mathbb{R}^{n\times \rho}}}{\int^{'}}
e^{i<x,\lambda^{'}_{1}+\cdots+\lambda^{'}_{\rho}>}
e^{-\mu\varepsilon
t\varepsilon^{-\frac{\alpha}{\alpha+\gamma}}|\lambda^{'}_{1}+\cdots+\lambda^{'}_{\rho}|^{\alpha}
(1+\varepsilon^{-\frac{2}{\alpha+\gamma}}|\lambda^{'}_{1}+\cdots+\lambda^{'}_{\rho}|^{2})^{\frac{\gamma}{2}}
}
\prod_{\sigma=1}^{\rho}{\sqrt{f(\varepsilon^{\chi}\lambda^{'}_{\sigma})}}
W(d\lambda^{'}_{\sigma}),
\end{align}
where we have used the self-similar property for Gaussian random
measure on $\mathbb{R}^{n}$ in the last equality.
Therefore, by the
orthogonal property for the Gaussian white noise,
we can get
\begin{align}\notag
&(C_{\rho}(h))^{-2}N_{1}(\varepsilon)\textup{Cov}(I_{\rho}^{\varepsilon}(t,x)I_{\rho}^{\varepsilon}(t^{'},x^{'}))
\\\notag
=&
\varepsilon^{\chi\rho n}\underset{{\mathbb{R}^{n\times \rho}}}{\int}
e^{i<x-x^{'},\lambda^{'}_{1}+\cdots+\lambda^{'}_{\rho}>}e^{ -\mu\varepsilon
(t+t^{'})\varepsilon^{-\frac{\alpha}{\alpha+\gamma}}|\lambda_{1}+\cdots+\lambda_{\rho}|^{\alpha}
(1+\varepsilon^{-\frac{2}{\alpha+\gamma}}|\lambda_{1}+\cdots+\lambda_{\rho}|^{2})^{\frac{\gamma}{2}}
}
\prod_{\sigma=1}^{\rho}{f(\varepsilon^{\chi}\lambda_{\sigma})}
d\lambda_{\sigma}\\\notag
=& \varepsilon^{\chi\rho n}
\underset{{\mathbb{R}^{n}}}{\int}
e^{i<x-x^{'},\tau_{1}>}
 e^{ -\mu\varepsilon
(t+t^{'})\varepsilon^{-\frac{\alpha}{\alpha+\gamma}}|\tau_{1}|^{\alpha}
(1+\varepsilon^{-\frac{2}{\alpha+\gamma}}|\tau_{1}|^{2})^{\frac{\gamma}{2}}
}
\frac{f^{*\rho}(\varepsilon^{\chi}\tau_{1})}{\varepsilon^{\chi(\rho-1)n}}d\tau_{1}\\\label{proof
micro single 2} =&
\varepsilon^{n\chi}\underset{{\mathbb{R}^{n}}}{\int}
e^{i<x-x^{'},\tau>} e^{ -\mu
(t+t^{'})|\tau|^{\alpha}
(\varepsilon^{\frac{2}{\alpha+\gamma}}+|\tau|^{2})^{\frac{\gamma}{2}}
} f^{*\rho}(\varepsilon^{\chi}\tau)d\tau,
\end{align}
where $f^{*\rho}(\cdot)$ is defined in (\ref{convolutionspectral}).\\
%%%%%%%%%%%%%%%%%%%%%%%% 0<\rho\alpha^{(j)}_{k}<n  %%%%%%%%%
(i) For $\rho\in \mathbb{N}$ with $m\kappa\leq\rho\kappa<n$
and any $\delta>0$, by (\ref{proof micro
single 2}) and (\ref{tauberian2}),
\begin{equation}\label{decom cov micro 1}
\textup{Cov}(I_{\rho}^{\varepsilon}(t,x)I_{\rho}^{\varepsilon}(t^{'},x^{'}))
=(C_{\rho}(h))^{2}(A_{1}(\varepsilon)+A_{2}(\varepsilon)),
\end{equation}
with
%\begin{align}\notag
%&(C_{\rho}(h))^{-2}\textup{Cov}(I_{\rho}^{\varepsilon}(t,x)I_{\rho}^{\varepsilon}(t^{'},x^{'}))
%\\\notag
%=&
%(N_{1}(\varepsilon))^{-1}\varepsilon^{n\chi}(1+o(1)) K(n,\rho\kappa)
%\underset{{\mathbb{R}^{n}}}{\int} e^{i<x-x^{'},\tau>}e^{ -\mu
%(t+t^{'})|\tau|^{\alpha}
%(\varepsilon^{\frac{2}{\alpha+\gamma}}+|\tau|^{2})^{\frac{\gamma}{2}}
%}\frac {L^{\rho}(|\varepsilon^{\chi}\tau|^{-1})} {|\varepsilon^{\chi}\tau|^{n-\rho\kappa}} d\tau
%\\\label{micro proof marker1}
%=&(1+o(1))(N_{1}(\varepsilon))^{-1} \varepsilon^{\chi\rho\kappa}L^{\rho}(\varepsilon^{-\chi})
%K(n,\rho\kappa)\underset{{\mathbb{R}^{n}}}{\int} e^{i<x-x^{'},\tau>}
%\frac{e^{ -\mu
%(t+t^{'})
%|\tau|^{\alpha+\gamma}}}{|\tau|^{n-\rho \kappa}}d\tau,
%\end{align}
\begin{align*}
|A_{1}(\varepsilon)|&=(N_{1}(\varepsilon))^{-1}\varepsilon^{n\chi}
|
\int_{{|\varepsilon^{\chi}\tau|>\delta}}
e^{i<x-x^{'},\tau>} e^{ -\mu
(t+t^{'})|\tau|^{\alpha}
(\varepsilon^{\frac{2}{\alpha+\gamma}}+|\tau|^{2})^{\frac{\gamma}{2}}
} f^{*\rho}(\varepsilon^{\chi}\tau)d\tau
|\\
&\leq
(N_{1}(\varepsilon))^{-1}\varepsilon^{n\chi}
\ \textup{sup}\{f^{*\rho}(\widetilde{\lambda})|\ |\widetilde{\lambda}|>\delta\}
\int_{{|\varepsilon^{\chi}\tau|>\delta}}
e^{-\mu (t+t^{'})|\tau|^{\alpha+\gamma}}d\tau\\
&\leq
(\varepsilon^{\chi m\kappa}L^{m}(\varepsilon^{-\chi}))^{-1}
\varepsilon^{n\chi}
\ \textup{sup}\{f^{*\rho}(\widetilde{\lambda})|\ |\widetilde{\lambda}|>\delta\}
\int_{\varepsilon^{-\chi}\delta}^{\infty}
e^{-\mu (t+t^{'})r^{\alpha+\gamma}}r^{n-1}dr\\
&\rightarrow 0,\ \textup{as}\ \varepsilon\rightarrow 0,\ \ (m\kappa\leq\rho\kappa<n)
\end{align*}
and, by choosing $\delta$ small enough and (\ref{tauberian2}),
\begin{align*}
&A_{2}(\varepsilon)\\
&=(N_{1}(\varepsilon))^{-1}\varepsilon^{n\chi}\int_{{|\varepsilon^{\chi}\tau|\leq\delta}}
e^{i<x-x^{'},\tau>} e^{ -\mu
(t+t^{'})|\tau|^{\alpha}
(\varepsilon^{\frac{2}{\alpha+\gamma}}+|\tau|^{2})^{\frac{\gamma}{2}}
} f^{*\rho}(\varepsilon^{\chi}\tau)d\tau\\
&
=
(N_{1}(\varepsilon))^{-1}\varepsilon^{n\chi}
\underset{|\varepsilon^{\chi}\tau|\leq\delta}{\int}
e^{i<x-x^{'},\tau>} e^{ -\mu
(t+t^{'})|\tau|^{\alpha}
(\varepsilon^{\frac{2}{\alpha+\gamma}}+|\tau|^{2})^{\frac{\gamma}{2}}
}
(1+o(1))
K(n,\rho\kappa)\frac{L^{\rho}(|\varepsilon^{\chi}\tau|^{-1})}
{|\varepsilon^{\chi}\tau|^{n-\rho\kappa}}d\tau
\\
&\ \sim
(N_{1}(\varepsilon))^{-1}\varepsilon^{\chi\rho\kappa}
L^{\rho}(\varepsilon^{-\chi})
K(n,\rho\kappa)
\underset{\mathbb{R}^{n}}{\int}
\frac
{
e^{i<x-x^{'},\tau>-\mu (t+t^{'})|\tau|^{\alpha+\gamma}}
}
{|\tau|^{n-\rho\kappa}}
d\tau,\ \textup{as}\ \varepsilon\rightarrow 0,
\end{align*}
where the asymptotic equivalence is guaranteed by the uniform
convergence theorem for the slowly varying function (see, for
example, Leonenko \cite[Section 1.4]{L}). We remark that the $f$,
defined as a spectral density function, is bounded outside of
zero(while the singularity at zero is from the LRD assumption, as
employed in subsection 3.2); therefore its $\rho$-th convolution
remains to have, at most, the singularity only at zero.
%
%where the argument about $``(1+o(1))"$ is guaranteed by the same method stated in
%
\\
The conclusion of (i): Apart from the term $\textup{Cov}
(I_{m}^{\varepsilon}(t,x)I_{m}^{\varepsilon}(t^{'},x^{'}))$,
\begin{align}\label{Cov decom micro 0}
\underset{\varepsilon\rightarrow 0}{\textup{lim}}
\underset{\rho: m<\rho<n/\kappa}{\sum}
\textup{Cov}
(I_{\rho}^{\varepsilon}(t,x)I_{\rho}^{\varepsilon}(t^{'},x^{'}))
=
 0
\end{align}
where we have used the fact that $\{l\in \mathbb{N}|m\kappa\leq\l\kappa<n\}$ is a finite set
and
on this set
$\underset{\varepsilon\rightarrow0}{\textup{lim}}(N_{1}(\varepsilon))^{-1} \varepsilon^{\chi\rho\kappa}L^{\rho}(\varepsilon^{-\chi})=0$
 except for $\rho=m$.\\
%
%%%%%%%%%%%%%%%%%%%%%%%%%%% \rho\\kappa^{(j)}_{0}>n%%%%%%%%%%%%%%%%%%%%%%%%%%%%%
(ii) For $\rho\in \mathbb{N}$ with $\rho\kappa>n$
, by (\ref{proof micro
single 2}), (\ref{L1spectral}) and
%\begin{align*}
%\textup{Cov}
%(I_{\rho}^{\varepsilon}(t,x)I_{\rho}^{\varepsilon}(t^{'},x^{'})) =&
%(1+o(1))\frac{\varepsilon^{n\chi}}{N_{1}(\varepsilon)}f^{*\rho}(0)(C_{\rho}(h))^{2}
%\underset{{\mathbb{R}^{n}}}{\int}e^{i<x-x^{'},\tau> -\mu
%(t+t^{'})
%|\tau|^{\alpha+\gamma}} d\tau.
%\end{align*}
%%%%%%%%%%%%%%%%%%%%%%%%%%%%%%%%%%%%%%%%%%%%%%%%%%%%%%
$\overset{\infty}{\underset{\rho=m}{\sum}}{(C_{\rho}(h))^{2}}\leq\parallel
h\parallel_{2}^{2}<\infty$,
\begin{align}\label{Cov decom micro 1}
&\underset{\varepsilon\rightarrow 0}{\textup{lim}}
\underset{\rho:\rho\kappa>n}{\sum}
\textup{Cov}
(I_{\rho}^{\varepsilon}(t,x)I_{\rho}^{\varepsilon}(t^{'},x^{'}))\\\notag
=
&\underset{\varepsilon\rightarrow 0}{\textup{lim}}\
\varepsilon^{n\chi}(N_{1}(\varepsilon))^{-1}
\underset{\rho:\rho\kappa>n}{\sum}(C_{\rho}(h))^{2}
\underset{\mathbb{R}^{n}}{\int}
e^{i<x-x^{'},\tau>-\mu (t+t^{'})|\tau|^{\alpha+\gamma}}
f^{*\rho}(0)
d\tau\\\notag
\leq&
\underset{\varepsilon\rightarrow 0}{\textup{lim}}\
\varepsilon^{n\chi}(N_{1}(\varepsilon))^{-1}
M\underset{\rho:\rho\kappa>n}{\sum}(C_{\rho}(h))^{2}
=0,
\end{align}
since by (\ref{L1spectral}) $f^{*\rho}(0)$ is bounded by $f^{*\widetilde{\rho}}(0)$
with $\widetilde{\rho}=\textup{inf}\{l\in\mathbb{N}|\ l\kappa>n\}$ so we set
$M:=f^{*\widetilde{\rho}}(0)\int_{\mathbb{R}^{n}}
e^{-\mu (t+t^{'})|\tau|^{\alpha+\gamma}}
d\tau.$
\\
Finally, from the expansion (\ref{expansion single micro}) for the random field
$s^{\varepsilon}(t,x)$
and combining the observations
(\ref{Cov decom micro 0})
and
(\ref{Cov decom micro 1})
we know that only the component $I_{m}^{\varepsilon}(t,x)$ in (\ref{expansion single micro})
do contribute to the covariance function of the random field $s^{\varepsilon}(t,x)$,
that is,
\begin{align*}
\underset{\varepsilon\rightarrow
0}{\textup{lim}}\ \textup{Cov}(s^{\varepsilon}(t,x)s^{\varepsilon}(t^{'},x^{'}))
=
(C_{m}(h))^{2}K(n,m\kappa)\underset{{\mathbb{R}^{n}}}{\int} e^{i<x-x^{'},\tau>}
 \frac{e^{ -\mu
(t+t^{'})
|\tau|^{\alpha+\gamma}}}{|\tau|^{n-m\kappa}}d\tau.
\ \ \ \ \ \ \ \ \ \ \ \ \ \ \  \raggedleft\qed\end{align*}
\\
%%%%%%%%%%%%%%%%%%%     %%%%%%%%%%%%%%%%%%%%%%%%        %%%%%%%%%%%%%%%%%%%%%%         %%%%%%%%%%
%%%%%%%%%%%%%%%%%%%%%%%%%%%%%%% Proof Prop 4 (2)  %%%%%%%%%%%%%%%%%%%%%%%%%%%%%%%%%%%%%%%%%%%%%%
{\bf (2)} From the above discussion, we may apply Chebyshev
inequality to obtain that:
\begin{align*}
\sum_{\rho=m+1}^{\infty}{I^{\varepsilon}_{\rho}(t,x)}\overset{P}{\longrightarrow}
0.
\end{align*}
Therefore, in view of Slutsky argument, we suffice to focus our
attention on the term $I^{\varepsilon}_{m}(t,x)$. In the following
we will prove $I^{\varepsilon}_{m}(t,x)$ converges in distribution
sense to $s_{m}(t,x)$, which is defined in (\ref{limiting
field(micro)}), for each fixed $(t,x)\in\mathbb{R}_{+}\times
\mathbb{R}^{n}$. By the definition of $N_{1}(\varepsilon)$ and
replacing the letter $\rho$ by $m$ in (\ref{proof micro single
1}), we can rewrite (\ref{proof micro single 1}) as follows
\begin{align}\label{proof micro single 5}
I^{\varepsilon}_{m}(t,x)
\overset{d}{=}
\frac{C_{m}(h)}{\sqrt{m!}}\underset{{\mathbb{R}^{n\times
m}}}{\int^{'}} e^{i<x,\lambda^{'}_{1}+\cdots+\lambda^{'}_{m}>}
M_{\varepsilon}(\lambda)
\prod_{\sigma=1}^{m}{W(d\lambda^{'}_{\sigma})},
\end{align}
with
\begin{align*}
M_{\varepsilon}(\lambda):= \varepsilon^{\frac{\chi m(n-\kappa)}{2}}
L^{-\frac{m}{2}}(\varepsilon^{-\chi}) e^{-\mu
t(\varepsilon^{\frac{2}{\alpha+\gamma}}+|\lambda_{1}+\cdots+\lambda_{m}|^{2})^{\frac{\gamma}{2}}
|\lambda_{1}+\cdots+\lambda_{m}|^{\alpha}}
\prod_{\sigma=1}^{m}{\sqrt{f(\varepsilon^{\chi}\lambda_{\sigma})}},
\end{align*}
which, when $\varepsilon\rightarrow 0$, satisfies
\begin{align}\label{proof micro single 6}
\underset{\varepsilon\rightarrow 0}{\textup{lim}}
M_{\varepsilon}(\lambda) \overset{(\ref{tauberian})}{=}
(K(n,\kappa))^{\frac{m}{2}}\frac{e^{-\mu t
|\lambda_{1}+\cdots\lambda_{m}|^{\alpha+\gamma}}}{(|\lambda_{1}|\cdots|\lambda_{m}|)^{\frac{n-\kappa}{2}}}.
\end{align}
%%%%%%%%%%%%%%%%%%%%%%%%%%%%%%%%%%%
Now, applying the isometric property of the multiple Wiener integrals to
the difference of (\ref{proof micro single 5}) and (\ref{limiting
field(micro)}), we have
\begin{align}\notag
E|I^{\varepsilon}_{m}(t,x) -s_{m}(t,x)|^{2}
=&\underset{\varepsilon\rightarrow
0}{\textup{lim}}(C_{\rho}(h))^{2}\underset{\mathbb{R}^{n\times
m}}{\int} |M_{\varepsilon}(\lambda) -
(K(n,\kappa))^{\frac{m}{2}}\frac{e^{-\mu t
|\lambda_{1}+\cdots\lambda_{m}|^{\alpha+\gamma}}}{(|\lambda_{1}|\cdots|\lambda_{m}|)^{\frac{n-\kappa}{2}}}|^{2}\prod_{\sigma=1}^{m}{d\lambda_{\sigma}}
\end{align}
$\longrightarrow 0,\ \textup{as}\ \varepsilon \rightarrow 0$, by
(\ref{proof micro single 6}) and the assumption $f(\lambda)$ is
decreasing at infinity in Condition C and
\begin{align*}
&\underset{\mathbb{R}^{n\times m}}{\int} \frac{e^{-2\mu t
|\lambda_{1}+\cdots\lambda_{m}|^{\alpha+\gamma}}}{(|\lambda_{1}|\cdots|\lambda_{m}|)^{n-\kappa}}\prod_{\sigma=1}^{m}{d\lambda_{\sigma}}
= r(n,m,\kappa)
\underset{\mathbb{R}^{n}}{\int}\frac{e^{-2\mu t
|\lambda|^{\alpha+\gamma}}}{(|\lambda|)^{n-m\kappa}}d\lambda<\infty,\
\textup{for}\ m\kappa<n,
\end{align*}
where the constant $r(n,m,\kappa)$ is generated by the Riesz
potential.
%For the technique for checking the legality of changing
%the order of limit and integration can refer to Leonenko \cite[p.
%252, 257]{L}.
Finally, the assertion (2) of Theorem 1 is followed from Slutsky
and Cramer-Wold arguments.
\qed\\
%
%
%%%%%%%%%%%%%%%%%%%%%%%%%%%%%%%%%%%%%%%%%%%%%%%%%%%%%%%%%%%%%
%%%%%%%%%%%%%%%%%%%%%%%%%%%%% Proof of Case 1. (micro, system )%%%%%%%%%%%%%%%%%%%%%%%%%%%%%%%%%%%%%%
                                 %%%%%%%%%%%%%%%%%%%%%%%%%%%%%%%%%%%%%%%%%
{\bf Proof of Theorem \ref{micro system theorem} for the case (1): $m_{2}\alpha_{2}>m_{1}\alpha_{1}.$}\\
From the solution form (\ref{solution}), we have
\begin{align}\label{proof micro system 1}
\left(\begin{array}{cc}u(t,x;u_{0}(\cdot))\\v(t,x;v_{0}(\cdot))\end{array}\right)-Q(t;d_{1},d_{2})\left(\begin{array}{cc}C_{0}^{(1)}\\C_{0}^{(2)}\end{array}\right)
=Q(t;d_{1},d_{2})\left[\left(\begin{array}{cccc}U(t,x)\\V(t,x)\end{array}\right)-\left(\begin{array}{cccc}C_{0}^{(1)}\\C_{0}^{(2)}\end{array}\right)\right]
,
\end{align}
where $Q(t;d_{1},d_{2})$, $U(t,x)$ and $V(t,x)$  are defined in (\ref{Q matrix}) and (\ref{homogeneous case}).\\
By (\ref{proof micro system 1}),
\begin{align}\notag
[\varepsilon^{m_{1}\kappa_{1}\chi}L^{m_{1}}(\varepsilon^{-\chi})]^{-\frac{1}{2}}
\left\{\left(\begin{array}{cccc} u(\varepsilon
t,\varepsilon^{\frac{1}{\alpha+\gamma}}x;u_{0}(\varepsilon^{-\frac{1}{\alpha+\gamma}-\chi}\cdot))\\v(\varepsilon
t,\varepsilon^{\frac{1}{\alpha+\gamma}}x;v_{0}(\varepsilon^{-\frac{1}{\alpha+\gamma}-\chi}\cdot))\end{array}\right)-Q(\varepsilon
t;d_{1},d_{2})\left(\begin{array}{cccc}C^{(1)}_{0}\\C^{(2)}_{0}
\end{array}\right)\right\}
\end{align}
\begin{align}
&=Q(\varepsilon
t;d_{1},d_{2})[\varepsilon^{m_{1}\kappa_{1}\chi}L^{m_{1}}(\varepsilon^{-\chi})]^{-\frac{1}{2}}
\left(\begin{array}{cc}U(\varepsilon
t,\varepsilon^{\frac{1}{\alpha+\gamma}}x;u_{0}(\varepsilon^{-\frac{1}{\alpha+\gamma}-\chi}\cdot))-C_{0}^{(1)}
\\V(\varepsilon
t,\varepsilon^{\frac{1}{\alpha+\gamma}}x;u_{0}(\varepsilon^{-\frac{1}{\alpha+\gamma}-\chi}\cdot))-C_{0}^{(2)}
\end{array}\right)\\
&:=Q(\varepsilon
t;d_{1},d_{2})[\varepsilon^{m_{1}\kappa_{1}\chi}L^{m_{1}}(\varepsilon^{-\chi})]^{-\frac{1}{2}}
\left(\begin{array}{cc} U_{\varepsilon}(t,x) \\ V_{\varepsilon}(t,x)
\end{array}\right).
\end{align}
Firstly, by Theorem \ref{Micro Scaling Theorem} (2), we have
\begin{align}\label{proof micro system 2}
U_{\varepsilon}(t,x)\Rightarrow \widetilde{X}^{(1)}_{m_{1}}(t,x),
\end{align}
where $\widetilde{X}^{(1)}_{m_{1}}$ is defined in (\ref{limit field micro case 1}).
\\
Secondly, by Theorem \ref{Micro Scaling Theorem} (1), we can obtain
\begin{align}\label{proof micro system 3}
V_{\varepsilon}(t,x)\overset{P}{\longrightarrow} 0
\end{align}
since we can apply Chebyshev inequality to observe that for any
$c>0$, as $\varepsilon\rightarrow 0$,
\begin{align}\notag
P(|V_{\varepsilon}(t,x)|>c)\leq c^{-2}
\textup{Var}(V_{\varepsilon}(t,x))\asymp
c^{-2}[\varepsilon^{-m_{1}\kappa_{1}\chi}L^{-m_{1}}(\varepsilon^{-\chi})]\cdot
[\varepsilon^{-m_{2}\kappa_{2}\chi}L^{-m_{2}}(\varepsilon^{-\chi})]
\rightarrow 0.
\end{align}
Meanwhile, since
\begin{align}\label{proof micro system 4}
\underset{\varepsilon\rightarrow 0}{\textup{lim}}\
Q_{\varepsilon}(t):= \underset{\varepsilon\rightarrow
0}{\textup{lim}}\ Q(\varepsilon
t;d_{1},d_{2})=P\left(\begin{array}{cccc}1 & 0
\\0 &
1\end{array}\right)P^{-1}=I_{2\times 2}.
\end{align}
\noindent Therefore, we may apply Lemma 1
 to those $U_{\varepsilon}(t,x),\ V_{\varepsilon}(t,x)$ and
$Q_{\varepsilon}(t)$ on the above to obtain that
\begin{align}\notag
&[\varepsilon^{m_{1}\kappa_{1}\chi}L^{m_{1}}(\varepsilon^{-\chi})]^{-\frac{1}{2}}
\left\{\left(\begin{array}{cccc} u(\varepsilon
t,\varepsilon^{\frac{1}{\alpha+\gamma}}x;u_{0}(\varepsilon^{-\frac{1}{\alpha+\gamma}-\chi}\cdot))\\v(\varepsilon
t,\varepsilon^{\frac{1}{\alpha+\gamma}}x;v_{0}(\varepsilon^{-\frac{1}{\alpha+\gamma}-\chi}\cdot))\end{array}\right)-Q(\varepsilon
t;d_{1},d_{2})\left(\begin{array}{cccc}C^{(1)}_{0}\\C^{(2)}_{0}
\end{array}\right)\right\}\\\notag
&\Rightarrow I_{2\times
2}\left(\begin{array}{cc}\widetilde{X}_{m_{1}}^{(1)}(t,x) \\
0\end{array}\right),\ \ t>0,\ x\in \mathbb{R}^{n}.
\hspace{7.9cm}\qed
\end{align}
                         %%%%%%%%%%%%%%%%%%%%%%%%%%%%%%%%
%%%%%%%%%%%%%%%%%%%%%%%%% Proof of Case 2. of system micro%%%%%%%%%%%%%%%%%%%%%%%%%%%
                     %%%%%%%%%%%%%%%%%%%%%%%%%%%%%%%%%%%%%%
{\bf Proof of Theorem \ref{micro system theorem} for the case (2): $m_{1}\alpha_{1}>m_{2}\alpha_{2}.$}\\
\noindent The proof is proceeded as the case (1),
 yet under the new assumption and the different renormalization
$[\varepsilon^{m_{2}\kappa_{2}\chi}L^{m_{2}}(\varepsilon^{-\chi})]^{-\frac{1}{2}}.$
Now (\ref{proof micro system 2}) becomes as
\begin{align*}
U_{\varepsilon}(t,x):=[\varepsilon^{m_{2}\kappa_{2}\chi}L^{m_{2}}(\varepsilon^{-\chi})]^{-\frac{1}{2}}
\Big\{ U(\varepsilon
t,\varepsilon^{\frac{1}{\alpha+\gamma}}x;u_{0}(\varepsilon^{-\frac{1}{\alpha+\gamma}-\chi}\cdot))-C_{0}^{(1)}\Big\}
\overset{P}{\longrightarrow} 0,
\end{align*}
and (\ref{proof micro system 3}) becomes as
\begin{align*}
V_{\varepsilon}(t,x):=
[\varepsilon^{m_{2}\kappa_{2}\chi}L^{m_{2}}(\varepsilon^{-\chi})]^{-\frac{1}{2}}
\Big\{ V(\varepsilon
t,\varepsilon^{\frac{1}{\alpha+\gamma}}x;v_{0}(\varepsilon^{-\frac{1}{\alpha+\gamma}-\chi}\cdot))-C_{0}^{(2)}\Big\}
\Rightarrow \widetilde{X}^{(2)}_{m_{2}}(t,x),
\end{align*}
where $\widetilde{X}^{(2)}_{m_{2}}$ is defined in (\ref{limit field micro case 2}).\\
While (\ref{proof micro system 4}) is kept unchanged. Therefore, we
again apply Lemma 1 to get
\begin{align}\notag
&[\varepsilon^{m_{2}\kappa_{2}\chi}L^{m_{2}}(\varepsilon^{-\chi})]^{-\frac{1}{2}}
\left\{\left(\begin{array}{cccc} u(\varepsilon
t,\varepsilon^{\frac{1}{\alpha+\gamma}}x;u_{0}(\varepsilon^{-\frac{1}{\alpha+\gamma}-\chi}\cdot))\\v(\varepsilon
t,\varepsilon^{\frac{1}{\alpha+\gamma}}x;v_{0}(\varepsilon^{-\frac{1}{\alpha+\gamma}-\chi}\cdot))\end{array}\right)-Q(\varepsilon
t;d_{1},d_{2})\left(\begin{array}{cccc}C^{(1)}_{0}\\C^{(2)}_{0}
\end{array}\right)\right\}\\\notag
&\Rightarrow I_{2\times 2}\left(\begin{array}{cc}0 \\
\widetilde{X}_{m_{2}}^{(2)}(t,x)\end{array}\right),\ \ t>0,\ x\in
\mathbb{R}^{n}.
\hspace{8.1cm}\qed\end{align}\\
                                  %%%%%%%%%%%%%%%%%%%%%%%%%%%%%%%%%%%%%%%%
%%%%%%%%%%%%%%%%%%%%%%%%%%% Proof of Case 3. of system  micro  %%%%%%%%%%%%%%%%%%%%%%%%
                               %%%%%%%%%%%%%%%%%%%%%%%%%%%%%%%%%%%%%%%
{\bf Proof of Theorem \ref{micro system theorem} for the case (3): $m_{1}=m_{2}=m,\ \alpha_{1}=\alpha_{2}=\alpha$.}\\
By  Theorem \ref{Micro Scaling Theorem} (2), we have
\begin{align*}
&U_{\varepsilon}(t,x):=[\varepsilon^{m\kappa\chi}L^{m}(\varepsilon^{-\chi})]^{-\frac{1}{2}}
\Big\{ U(\varepsilon
t,\varepsilon^{\frac{1}{\alpha+\gamma}}x;u_{0}(\varepsilon^{-\frac{1}{\alpha+\gamma}-\chi}\cdot))-C_{0}^{(1)}\Big\}\Rightarrow
\widetilde{X}^{(1)}_{m}(t,x),\\
&V_{\varepsilon}(t,x):=
[\varepsilon^{m\kappa\chi}L^{m}(\varepsilon^{-\chi})]^{-\frac{1}{2}}
\Big\{ V(\varepsilon
t,\varepsilon^{\frac{1}{\alpha+\gamma}}x;v_{0}(\varepsilon^{-\frac{1}{\alpha+\gamma}-\chi}\cdot))-C_{0}^{(2)}\Big\}\Rightarrow
\widetilde{X}^{(2)}_{m}(t,x),
\end{align*}
where $\widetilde{X}^{(j)}_{m},\ j\in\{1,2\}$, are defined in
(\ref{limit field micro case 1}) and (\ref{limit field micro case 2}) with $m_{1}=m_{2}=m$.
\\
Because in this case the equality
$\underset{\varepsilon\rightarrow0}{\textup{lim}} Q_{\varepsilon}(t)=I$ is still unchange,
in the same way, we obtained
\begin{align}\notag
&[\varepsilon^{m\kappa\chi}L^{m}(\varepsilon^{-\chi})]^{-\frac{1}{2}}
\left\{\left(\begin{array}{cccc} u(\varepsilon
t,\varepsilon^{\frac{1}{\alpha+\gamma}}x;u_{0}(\varepsilon^{-\frac{1}{\alpha+\gamma}-\chi}\cdot))\\v(\varepsilon
t,\varepsilon^{\frac{1}{\alpha+\gamma}}x;v_{0}(\varepsilon^{-\frac{1}{\alpha+\gamma}-\chi}\cdot))\end{array}\right)-Q(\varepsilon
t;d_{1},d_{2})\left(\begin{array}{cccc}C^{(1)}_{0}\\C^{(2)}_{0}
\end{array}\right)\right\}
\\\notag
&\Rightarrow I_{2\times 2}
\left(\begin{array}{cc}\widetilde{X}_{m}^{(1)}(t,x) \\
\widetilde{X}_{m}^{(2)}(t,x)\end{array}\right) =
\left(\begin{array}{cc} \widetilde{X}_{m}^{(1)}(t,x) \\
\widetilde{X}_{m}^{(2)}(t,x)\end{array}\right),\ \ t>0,\ x\in
\mathbb{R}^{n}.\hspace{4.8cm}\qed
\end{align}
%%%%%%%%%%%%%%%%%%%%%%%%%%%%%%%%%%%%%%%%%%%%%%%%%%%%%%%%%%%%%%%%%
%%%%%%%%%%%%%%%%%%%%%%%%%%%%%% %%%%%%%%%%%%%%%%%%%%%%%%%% %%%%%%%%%%%%%%%%%%%%%%%%%%% %%%%%%%%%%%%%%%%%%%%%%%%%%%
%%%%%%%%%%%%%%%% %%%%%%%%%%%%%%%%%%%%%%%% Proof of Proposition 5 %%%%%%%%%%%%%%%%%%%%% %%%%%%%%%%%%%%%%%%%%%%%%%%
\noindent{\bf Proof of Proposition \ref{Cov matrix limit micro
system}.}
\begin{align}\notag
&E\left(\begin{array}{cc}Y^{***}_{1}(t,x)\\Y^{***}_{2}(t,x)\end{array}\right)
\left(\begin{array}{cc}Y^{***}_{1}(t^{'},x^{'}) &
Y^{***}_{2}(t^{'},x^{'})\end{array}\right)\\\notag
=&E
\left(\begin{array}{cc}\widetilde{X}^{(1)}_{m}(t,x)\\\widetilde{X}^{(2)}_{m}(t,x)\end{array}\right)
\left(\begin{array}{cc}\widetilde{X}^{(1)}_{m}(t^{'},x^{'}) &
\widetilde{X}^{(2)}_{m}(t^{'},x^{'})\end{array}\right)\\\label{Spectral calculus}
=&\left(\begin{array}{cccc}E\widetilde{X}^{(1)}_{m}(t,x)\widetilde{X}^{(1)}_{m}(t^{'},x^{'})&
0 \\ 0 & E\widetilde{X}^{(2)}_{m}(t,x)\widetilde{X}^{(2)}_{m}(t^{'},x^{'}).\end{array}\right)
\end{align}
Because the representation for the limiting fields $\widetilde{X}^{(1)}_{m}(t,x)$
and
$\widetilde{X}^{(2)}_{m}(t,x)$ is the same as the limiting field $s_{m}(t,x)$,
defined in
(\ref{limiting field(micro)}),
we can apply the result  (\ref{cov micro single thm}) to get
\begin{align*}
E\widetilde{X}^{(j)}_{m}(t,x)\widetilde{X}^{(j)}_{m}(t^{'},x^{'})=
(C^{(j)}_{m})^{2}K(n,\kappa m)\underset{{\mathbb{R}^{n}}}{\int} e^{i<x-x^{'},\tau>}
 \frac{e^{ -\mu
(t+t^{'})
|\tau|^{\alpha+\gamma}}}{|\tau|^{n-m\kappa}}d\tau.
\end{align*}
Therefore, the covariance structure (\ref{Spectral calculus})
is equal to
\[
\ \ \ \ \ \ \  \ \ \ \ \ \ \ \ \ \  \
\int_{\mathbb{R}^{n}}
e^{i<x-x^{'},\tau>}
 K(n,\kappa m)\frac{e^{ -\mu
(t+t^{'})
|\tau|^{\alpha+\gamma}}}{|\tau|^{n-m\kappa}}
\left(
\begin{array}{cccc}
(C^{(1)}_{m})^{2} & 0 \\ 0 & (C^{(2)}_{m})^{2}
\end{array}
\right) d\tau. \ \ \ \ \ \ \ \ \ \ \ \ \qed
\]
%%%%%%%%%%%%%%%%%%%%%%%%%%%%%%%%%%%%%%%%%%%%%%%%%%%%%%%%%%%%%%%%%%%%%%
%
%
%%%%%%%%%%%%% Proof of the results about Macro single equation. %%%%%%%%%
\noindent{\bf Proof of Proposition \ref{cite theorem (Anh)}}\\
{\bf (1)} Here, for simplification, we set
$G(t,x):=G(t,x;\alpha,\gamma)$ and
$N(\varepsilon):=\varepsilon^{\frac{m\kappa}{\alpha}}L^{m}(\varepsilon^{-\frac{1}{\alpha}})$
, then by the solution form (\ref{homogeneous case}) for the
differential equation $\frac{\partial}{\partial
t}s=-\mu(I-\Delta)^{\frac{\gamma}{2}} (-\Delta)^{\frac{\alpha}{2}}
s$ we have
\begin{align}\notag
(N(\varepsilon))^{-\frac{1}{2}}&
\Big\{
s(\frac{t}{\varepsilon},\frac{x}{\varepsilon^{\frac{1}{\alpha}}};h(\zeta(\cdot)))-C_{0}(h)
\Big\}
\\\notag
=&
(N(\varepsilon))^{-\frac{1}{2}}
\Big\{
\int_{\mathbb{R}^{n}}
G(\frac{t}{\varepsilon},\frac{x}{\varepsilon^{\frac{1}{\alpha}}}-y)
\Bigl[
C_{0}(h)
+\overset{\infty}{\underset{k=m}{\sum}}
C_{k}(h)\frac{H_{k}(\zeta(y))}{\sqrt{k!}}
\Bigr]
dy-C_{0}(h)
\Big\}\\\label{decomposition solution}
\overset{(\ref{integral of G})}{=}&
\overset{\infty}{\underset{k=m}{\sum}}
\frac{C_{k}(h)}{\sqrt{k!}}(N(\varepsilon))^{-\frac{1}{2}}
\int_{\mathbb{R}^{n}}
G(\frac{t}{\varepsilon},\frac{x}{\varepsilon^{\frac{1}{\alpha}}}-y)
H_{k}(\zeta(y))dy:=\overset{\infty}{\underset{k=m}{\sum}}s^{\varepsilon}_{k}(t,x).
\end{align}
From (\ref{expectionhermite}), the cross terms of the left hand side blow have zero covariance,
thus we have
\begin{align}\label{Cov decom 0}
\textup{Cov}(\overset{\infty}{\underset{k=m}{\sum}}s^{\varepsilon}_{k}(t,x),
\overset{\infty}{\underset{k=m}{\sum}}s^{\varepsilon}_{k}(t^{'},x^{'}))
=
\overset{\infty}{\underset{k=m}{\sum}}
\textup{Cov}(s^{\varepsilon}_{k}(t,x),s^{\varepsilon}_{k}(t^{'},x^{'})).
\end{align}
For each $k\in \{m,m+1,\ldots\}$, by (\ref{expectionhermite})
\begin{align}\notag
&\textup{Cov}(s^{\varepsilon}_{k}(t,x),s^{\varepsilon}_{k}(t^{'},x^{'}))
\\\notag
=&
(C_{k}(h))^{2}(N(\varepsilon))^{-1}
\underset{\mathbb{R}^{n}}{\int}\underset{\mathbb{R}^{n}}{\int}
G(\frac{t}{\varepsilon},\frac{x}{\varepsilon^{\frac{1}{\alpha}}}-y)
G(\frac{t^{'}}{\varepsilon},\frac{x^{'}}{\varepsilon^{\frac{1}{\alpha}}}-y^{'})
R^{k}(y-y^{'})dy\ dy^{'}\\\notag
=\ &
(C_{k}(h))^{2}(N(\varepsilon))^{-1}
\underset{\mathbb{R}^{n}}{\int}\underset{\mathbb{R}^{n}}{\int}
G(\frac{t}{\varepsilon},\frac{x}{\varepsilon^{\frac{1}{\alpha}}}-y)
G(\frac{t^{'}}{\varepsilon},\frac{x^{'}}{\varepsilon^{\frac{1}{\alpha}}}-y^{'})
\Big\{
\underset{\mathbb{R}^{n}}{\int}
e^{i<y-y^{'},\lambda>}f^{*k}(\lambda)d\lambda\Big\}dy\ dy^{'}
\\\notag
=\ &
(C_{k}(h))^{2}(N(\varepsilon))^{-1}
\underset{\mathbb{R}^{n}}{\int}
\Big\{
\underset{\mathbb{R}^{n}}{\int}
e^{i<y,\lambda>}G(\frac{t}{\varepsilon},\frac{x}{\varepsilon^{\frac{1}{\alpha}}}-y)
dy
\underset{\mathbb{R}^{n}}{\int}
e^{i<-y^{'},\lambda>}G(\frac{t^{'}}{\varepsilon},\frac{x^{'}}{\varepsilon^{\frac{1}{\alpha}}}-y^{'})
dy^{'}
\Big\}
f^{*k}(\lambda)d\lambda
\\\notag
\overset{(\ref{fractional Green function})}{=}&
(C_{k}(h))^{2}(N(\varepsilon))^{-1}
\underset{\mathbb{R}^{n}}{\int}
e^{i<\frac{x-x^{'}}{\varepsilon^{\frac{1}{\alpha}}},\lambda>}
e^{-\mu\frac{t+t^{'}}{\varepsilon}|\lambda|^{\alpha}(1+|\lambda|^{2})^{\frac{\gamma}{2}}}
f^{*k}(\lambda)d\lambda
\\\label{Cov decomposition}
=&
(C_{k}(h))^{2}(N(\varepsilon))^{-1}
\underset{\mathbb{R}^{n}}{\int}
e^{i<x-x^{'},\lambda>}
e^{-\mu (t+t^{'})|\lambda|^{\alpha}(1+|\varepsilon^{\frac{1}{\alpha}}\lambda|^{2})^{\frac{\gamma}{2}}}
\varepsilon^{\frac{n}{\alpha}}f^{*k}(\varepsilon^{\frac{1}{\alpha}}\lambda)d\lambda,
\end{align}
by rescaling $\lambda$ into $\varepsilon^{\frac{1}{\alpha}}\lambda$.\\
For $k:m\kappa\leq k\kappa<n$ and any $\delta>0$,
by (\ref{Cov decomposition})
\begin{align}\label{Cov decom 1}
\textup{Cov}(s^{\varepsilon}_{k}(t,x),s^{\varepsilon}_{k}(t^{'},x^{'}))
=(C_{k}(h))^{2}(A_{1}(\varepsilon)+A_{2}(\varepsilon)),
\end{align}
with\\[-1.2cm]
\begin{align*}
|A_{1}(\varepsilon)|&=(N(\varepsilon))^{-1}
|
\int_{{|\varepsilon^{\frac{1}{\alpha}}\lambda|>\delta}}
e^{i<x-x^{'},\lambda>}
e^{-\mu (t+t^{'})|\lambda|^{\alpha}(1+|\varepsilon^{\frac{1}{\alpha}}\lambda|^{2})^{\frac{\gamma}{2}}}
\varepsilon^{\frac{n}{\alpha}}f^{*k}(\varepsilon^{\frac{1}{\alpha}}\lambda)d\lambda
|\\
&\leq
(N(\varepsilon))^{-1}\varepsilon^{\frac{n}{\alpha}}
\ \textup{sup}\{f^{*k}(\widetilde{\lambda})|\ |\widetilde{\lambda}|>\delta\}
\int_{{|\varepsilon^{\frac{1}{\alpha}}\lambda|>\delta}}
e^{-\mu (t+t^{'})|\lambda|^{\alpha}}d\lambda\\
&\leq
C(\delta)\ \varepsilon^{\frac{n-m\kappa}{\alpha}}L^{m}(\varepsilon^{-\frac{1}{\alpha}})
\int_{\varepsilon^{-\frac{1}{\alpha}}\delta}^{\infty}
e^{-\mu (t+t^{'})r^{\alpha}}r^{n-1}dr
\rightarrow 0,\ \textup{as}\ \varepsilon\rightarrow 0,
\end{align*}
and by choosing $\delta$ small enough
\begin{align*}
A_{2}(\varepsilon)
&=(N(\varepsilon))^{-1}\varepsilon^{\frac{n}{\alpha}}
\underset{|\varepsilon^{\frac{1}{\alpha}}\lambda|\leq\delta}{\int}
e^{i<x-x^{'},\lambda>-\mu (t+t^{'})|\lambda|^{\alpha}(1+|\varepsilon^{\frac{1}{\alpha}}\lambda|^{2})^{\frac{\gamma}{2}}}
f^{*k}(\varepsilon^{\frac{1}{\alpha}}\lambda)d\lambda\\
\end{align*}
\begin{align*}
&
\overset{(\ref{tauberian2})}{=}
(N(\varepsilon))^{-1}\varepsilon^{\frac{n}{\alpha}}
\underset{|\varepsilon^{\frac{1}{\alpha}}\lambda|\leq\delta}{\int}
e^{i<x-x^{'},\lambda>-\mu (t+t^{'})|\lambda|^{\alpha}(1+|\varepsilon^{\frac{1}{\alpha}}\lambda|^{2})^{\frac{\gamma}{2}}}
(1+o(1))
K(n,k\kappa)\frac{L^{k}(|\varepsilon^{\frac{1}{\alpha}}\lambda|^{-1})}
{|\varepsilon^{\frac{1}{\alpha}}\lambda|^{n-k\kappa}}d\lambda
\\
&\sim
(N(\varepsilon))^{-1}\varepsilon^{\frac{k\kappa}{\alpha}}
L^{k}(\varepsilon^{-\frac{1}{\alpha}})
(1+\widetilde{\varepsilon})K(n,k\kappa)
\underset{\mathbb{R}^{n}}{\int}
\frac
{
e^{i<x-x^{'},\lambda>-\mu (t+t^{'})|\lambda|^{\alpha}}
}
{|\lambda|^{n-k\kappa}}
d\lambda,\ \textup{as}\ \varepsilon\rightarrow 0,
\end{align*}
where the asymptotic equivalence is guaranteed by the uniform
convergence theorem for the slowly varying function (see, for
example, \cite[Section 1.4]{L}).
So from the above discussions in (\ref{Cov decom 1}) we can conclude that
\begin{align}\label{remainder term pro3}
\underset{\varepsilon\rightarrow 0}{\textup{lim}}
\underset{k:m\kappa\leq k\kappa<n}{\sum}
\textup{Cov}(s^{\varepsilon}_{k}(t,x),s^{\varepsilon}_{k}(t^{'},x^{'}))
=(C_{m}(h))^{2}K(n,m\kappa)
\underset{\mathbb{R}^{n}}{\int}
\frac
{
e^{i<x-x^{'},\lambda>-\mu (t+t^{'})|\lambda|^{\alpha}}
}
{|\lambda|^{n-m\kappa}}
d\lambda,
\end{align}
since $\{l\in \mathbb{N}|m\kappa\leq l\kappa<n\}$ is a finite set and
$\underset{\varepsilon\rightarrow 0}{\textup{lim}}
(N(\varepsilon))^{-1}\varepsilon^{\frac{k\kappa}{\alpha}}L^{k}(\varepsilon^{-\frac{1}{\alpha}})=0$
for $k\in\{l\in \mathbb{N}|m\kappa\leq l\kappa<n\}$ except for the term $k=m$.
\\
For $k:k\kappa>n$, by (\ref{Cov decomposition}) and (\ref{L1spectral})
\begin{align}\label{Cov decom 2}
\underset{\varepsilon\rightarrow 0}{\textup{lim}}
\underset{k:k\kappa>n}{\sum}
&\textup{Cov}(s^{\varepsilon}_{k}(t,x),s^{\varepsilon}_{k}(t^{'},x^{'}))
=
\\\notag
&\underset{\varepsilon\rightarrow 0}{\textup{lim}}\varepsilon^{\frac{n}{\alpha}}(N(\varepsilon))^{-1}
(1+o(1))\underset{k:k\kappa>n}{\sum}(C_{k}(h))^{2}
\underset{\mathbb{R}^{n}}{\int}
e^{i<x-x^{'},\lambda>-\mu (t+t^{'})|\lambda|^{\alpha}}
f^{*k}(0)
d\lambda\\\notag
&\underset{\varepsilon\rightarrow 0}{\textup{lim}}\varepsilon^{\frac{n}{\alpha}}(N(\varepsilon))^{-1}
(1+o(1))M\underset{k:k\kappa>n}{\sum}(C_{k}(h))^{2}
=0,
\end{align}
since by (\ref{L1spectral}) $f^{*k}(0)$ is bounded by $f^{*\widetilde{k}}(0)$
with $\widetilde{k}=\textup{inf}\{l\in\mathbb{N}|\ l\kappa>n\}$ so we set
$M:=f^{*\widetilde{k}}(0)\int_{\mathbb{R}^{n}}
e^{-\mu (t+t^{'})|\lambda|^{\alpha}}
d\lambda.$
The proof of Proposition \ref{cite theorem (Anh)} (1) is completed by
combining (\ref{Cov decom 0}), (\ref{remainder term pro3}) and (\ref{Cov decom 2})
to obtain
\begin{align*}
\underset{\varepsilon\rightarrow 0}{\textup{lim}}\
\textup{Cov}(s^{\varepsilon}(t,x),s^{\varepsilon}(t^{'},x^{'}))=
(C_{m}(h))^{2}K(n,m\kappa)
\underset{\mathbb{R}^{n}}{\int}
\frac
{
e^{i<x-x^{'},\lambda>-\mu (t+t^{'})|\lambda|^{\alpha}}
}
{|\lambda|^{n-m\kappa}}
d\lambda.\hspace{2.2cm}\qed
\end{align*}
\noindent {\bf (2)} This is derived from (1) by the same way as
that in the proof of \cite[Theorems 2.2 and 2.3]{AnhHomo}, and
thus we omit it. \qed

%           %%%%%%%%%%%%%%%%%%%%%%%%%%%%%%%%%%%%%%%%%%%
%%%%%%%%%%% Proof of Case 1. of Theorem  Macro system %%%%%%%%%%%%%%%%%%%%
           %%%%%%%%%%%%%%%%%%%%%%%%%%%%%%%%%%%%%%%%%
\noindent {\bf Proof of Theorem \ref{macro system theorem} for the case (1):} $m_{2}\kappa_{2}>m_{1}\kappa_{1}\ \textup{and}\ d_{1}>d_{2}.$\\
From (\ref{proof micro system 1}), we have
\begin{align}\notag
[\varepsilon^{\frac{m_{1}\kappa_{1}}{\alpha}}L^{m_{1}}(\varepsilon^{-\frac{1}{\alpha}})]^{-\frac{1}{2}}e^{-d_{1}\frac{t}{\varepsilon}}
\left\{\left(\begin{array}{cccc}u(\frac{t}{\varepsilon},\frac{x}{\varepsilon^{\frac{1}{\alpha}}})
\\v(\frac{t}{\varepsilon},\frac{x}{\varepsilon^{\frac{1}{\alpha}}})\end{array}\right)-Q(\frac{t}{\varepsilon};d_{1},d_{2})\left(\begin{array}{cccc}C^{(1)}_{0}\\C^{(2)}_{0}
\end{array}\right)\right\}
\end{align}
\begin{align}\label{Macro system proof 1}
=e^{-d_{1}\frac{t}{\varepsilon}}Q(\frac{t}{\varepsilon};d_{1},d_{2})
[\varepsilon^{\frac{m_{1}\kappa_{1}}{\alpha}}L^{m_{1}}(\varepsilon^{-\frac{1}{\alpha}})]^{-\frac{1}{2}}
\left(\begin{array}{cccc}U(\frac{t}{\varepsilon},\frac{x}{\varepsilon^{\frac{1}{\alpha}}})-C_{0}^{(1)}
\\
V(\frac{t}{\varepsilon},\frac{x}{\varepsilon^{\frac{1}{\alpha}}})-C_{0}^{(2)}\end{array}\right).
\end{align}
Firstly, by Proposition \ref{cite theorem (Anh)} (2), we have
\begin{align}
U_{\varepsilon}(t,x):=[\varepsilon^{\frac{m_{1}\kappa_{1}}{\alpha}}L^{m_{1}}(\varepsilon^{-\frac{1}{\alpha}})]^{-\frac{1}{2}}
\Big\{U(\frac{t}{\varepsilon},\frac{x}{\varepsilon^{\frac{1}{\alpha}}})-C_{0}^{(1)}\Big\}\Longrightarrow
X^{(1)}_{m_{1}}(t,x),
\end{align}
where $X^{(1)}_{m_{1}}$ is defined in (\ref{macro thm 1}).
\\
Secondly, by Proposition \ref{cite theorem (Anh)} (1), we have
\begin{align}\label{Macro system case 1.1}
V_{\varepsilon}(t,x)
:=[\varepsilon^{\frac{m_{1}\kappa_{1}}{\alpha}}L^{m_{1}}(\varepsilon^{-\frac{1}{\alpha}})]^
{-\frac{1}{2}}\Big\{V(\frac{t}{\varepsilon},\frac{x}{\varepsilon^{\frac{1}{\alpha}}})-C_{0}^{(2)}\Big\}
\overset{P}{\longrightarrow} 0,
\end{align}
since we can apply  Chebyshev inequality to observe that, for any
$c>0$, as $\varepsilon\rightarrow 0$,
\begin{align}\notag
P(|V_{\varepsilon}(t,x)|>c)\leq c^{-2}
\textup{Var}(V_{\varepsilon}(t,x))\asymp
c^{-2}[\varepsilon^{-\frac{m_{1}\kappa_{1}}{\alpha}}L^{-m_{1}}(\varepsilon^{-\frac{1}{\alpha}})]\cdot
[\varepsilon^{\frac{m_{2}\kappa_{2}}{\alpha}}L^{m_{2}}(\varepsilon^{-\frac{1}{\alpha}})]
\rightarrow 0 .
\end{align}
Meanwhile, since $d_{1}>d_{2},$
\begin{align}\notag
\underset{\varepsilon\rightarrow 0}{\textup{lim}}\
Q_{\varepsilon}(t):&= \underset{\varepsilon\rightarrow
0}{\textup{lim}}\
e^{-d_{1}\frac{t}{\varepsilon}}Q(\frac{t}{\varepsilon};d_{1},d_{2})
=\underset{\varepsilon\rightarrow
0}{\textup{lim}}\
e^{-d_{1}\frac{t}{\varepsilon}}P\left(\begin{array}{cccc}e^{-d_{1}\frac{t}{\varepsilon}}
& 0 \\0 &
e^{-d_{2}\frac{t}{\varepsilon}}\end{array}\right)P^{-1}\\\label{limit matrix}
&=P\left(\begin{array}{cccc}1
& 0 \\0 &
0\end{array}\right)P^{-1}
=\left(\begin{array}{cccc}p_{11}p_{22} &
-p_{11}p_{12}\\p_{21}p_{22} & -p_{12}p_{21}\end{array}\right).
\end{align}
Therefore,  by  the independence between
$U_{\varepsilon}(t,x)$ and $V_{\varepsilon}(t,x)$  we may apply Lemma 1  to obtain
\begin{align}\notag
&[\varepsilon^{\frac{m_{1}\kappa_{1}}{\alpha}}L^{m_{1}}(\varepsilon^{-\frac{1}{\alpha}})]^{-\frac{1}{2}}e^{-d_{1}\frac{t}{\varepsilon}}
\left\{\left(\begin{array}{cc}u(\frac{t}{\varepsilon},\frac{x}{\varepsilon^{\frac{1}{\alpha}}})
\\v(\frac{t}{\varepsilon},\frac{x}{\varepsilon^{\frac{1}{\alpha}}})\end{array}\right)-Q(\frac{t}{\varepsilon};d_{1},d_{2})\left(\begin{array}{cccc}C^{(1)}_{0}\\C^{(2)}_{0}
\end{array}\right)\right\}\\\notag
&\Rightarrow\left(\begin{array}{cccc}p_{11}p_{22} &
-p_{11}p_{12}\\p_{21}p_{22} & -p_{12}p_{21}\end{array}\right)
\left(\begin{array}{cc}X_{m_{1}}^{(1)}(t,x) \\
0\end{array}\right)\\\notag
&=\left(\begin{array}{cccc}p_{11}p_{22}X^{(1)}_{m_{1}}(t,x)\\p_{21}p_{22}X^{(1)}_{m_{1}}(t,x)\end{array}\right),\
\ t>0,\ x\in \mathbb{R}^{n}.\hspace{7.5cm}\qed
\end{align}
%
%%%%%%%%%%%%%%%%%%%%%%%%%%%%%%%%%%%%%%% proof Case (2) Macro system $$$$$$$$$$$$$$$$$$$$$$$$$$$$$$$$$$$$
%
{\bf Proof of Theorem \ref{macro system theorem} for the case (2):} $m_{2}\kappa_{2}<m_{1}\kappa_{1}\ \textup{and}\ d_{1}>d_{2}.$\\
We use the same scheme as in the proof of the case (1) with the roles of
$U_{\varepsilon}(t,x)$ and $V_{\varepsilon}(t,x)$ being replaced as follows
\begin{align}\label{Macro system case 2.2}
U_{\varepsilon}(t,x)
:=[\varepsilon^{\frac{m_{2}\kappa_{2}}{\alpha}}L^{m_{2}}(\varepsilon^{-\frac{1}{\alpha}})]^
{-\frac{1}{2}}\Big\{U(\frac{t}{\varepsilon},\frac{x}{\varepsilon^{\frac{1}{\alpha}}})-C_{0}^{(1)}\Big\}
\overset{P}{\longrightarrow} 0.
\end{align}
and
\begin{align}\label{Macro system case 2.1}
V_{\varepsilon}(t,x):=[\varepsilon^{\frac{m_{2}\kappa_{2}}{\alpha}}L^{m_{2}}(\varepsilon^{-\frac{1}{\alpha}})]^{-\frac{1}{2}}
\Big\{V(\frac{t}{\varepsilon},\frac{x}{\varepsilon^{\frac{1}{\alpha}}})-C_{0}^{(2)}\Big\}\Longrightarrow
X^{(2)}_{m_{2}}(t,x),
\end{align}
where $X^{(2)}_{m_{2}}$ is defined in (\ref{macro thm 2}).
Additionally, in this case the limit of the matrix
$Q_{\varepsilon}(t)$ coincides with (\ref{limit matrix})
so
\begin{align}\notag
&[\varepsilon^{\frac{m_{2}\kappa_{2}}{\alpha}}L^{m_{2}}(\varepsilon^{-\frac{1}{\alpha}})]^{-\frac{1}{2}}e^{-d_{1}\frac{t}{\varepsilon}}
\left\{\left(\begin{array}{cc}u(\frac{t}{\varepsilon},\frac{x}{\varepsilon^{\frac{1}{\alpha}}})
\\v(\frac{t}{\varepsilon},\frac{x}{\varepsilon^{\frac{1}{\alpha}}})\end{array}\right)
-Q(\frac{t}{\varepsilon};d_{1},d_{2})\left(\begin{array}{cc}C^{(1)}_{0}\\C^{(2)}_{0}
\end{array}\right)\right\}\\\notag
&\Rightarrow\left(\begin{array}{cccc}p_{11}p_{22} &
-p_{11}p_{12}\\p_{21}p_{22} & -p_{12}p_{21}\end{array}\right)
\left(\begin{array}{cc}0 \\
X_{m_{2}}^{(2)}(t,x)\end{array}\right)\\\notag
&=\left(\begin{array}{cc}-p_{11}p_{12}X^{(2)}_{m_{2}}(t,x)\\-p_{12}p_{21}X^{(2)}_{m_{2}}(t,x)\end{array}\right),\
\ t>0,\ x\in \mathbb{R}^{n}.
\hspace{7.5cm}\qed\end{align}
\\
%%%%%%%%%%%%%%%%%%%%% proof Case (3) Macro system %%%%%%%%%%%%%%%%%%%%%%%%%%%%%%%%
%
{\bf Proof of Theorem \ref{macro system theorem} for the case (3):} $m_{1}=m_{2}=m$,
$\kappa_{1}=\kappa_{2}=\kappa$
and $d_{1}>d_{2}.$\\
By Proposition \ref{cite theorem (Anh)} (2), we have
\begin{align}\label{macro case3 U}
U_{\varepsilon}(t,x):=[\varepsilon^{\frac{m\kappa}{\alpha}}L^{m}(\varepsilon^{-\frac{1}{\alpha}})]^{-\frac{1}{2}}
\Big\{U(\frac{t}{\varepsilon},\frac{x}{\varepsilon^{\frac{1}{\alpha}}})-C_{0}^{(1)}\Big\}\Longrightarrow
X^{(1)}_{m}(t,x),
\end{align}
and
\begin{align}\label{macro case3 V}
V_{\varepsilon}(t,x):=[\varepsilon^{\frac{m\kappa}{\alpha}}L^{m}(\varepsilon^{-\frac{1}{\alpha}})]^{-\frac{1}{2}}
\Big\{V(\frac{t}{\varepsilon},\frac{x}{\varepsilon^{\frac{1}{\alpha}}})-C_{0}^{(2)}\Big\}\Longrightarrow
X^{(2)}_{m}(t,x),
\end{align}
where $X^{(1)}_{m}$ and $X^{(2)}_{m}$ is defined in (\ref{macro thm 1}) and (\ref{macro thm 2})
with $m_{1}=m_{2}=m$.
\\
Additionally, in this case the matrix $Q_{\varepsilon}(t)$ is also unchanged
so
by applying Lemma 1 to $U_{\varepsilon}(t,x)$,
$V_{\varepsilon}(t,x)$ and $Q_{\varepsilon}(t)$ which are given in
(\ref{macro case3 U}), (\ref{macro case3 V}) and (\ref{limit matrix}), respectively, we see that the finite
dimensional distributions of the rescaled random field
\begin{align}\notag
&[\varepsilon^{\frac{m\kappa}{\alpha}}L^{m}(\varepsilon^{-\frac{1}{\alpha}})]^{-\frac{1}{2}}e^{-d_{1}\frac{t}{\varepsilon}}
\left\{\left(\begin{array}{cc}u(\frac{t}{\varepsilon},\frac{x}{\varepsilon^{\frac{1}{\alpha}}})
\\v(\frac{t}{\varepsilon},\frac{x}{\varepsilon^{\frac{1}{\alpha}}})\end{array}\right)
-Q(\frac{t}{\varepsilon};d_{1},d_{2})\left(\begin{array}{cc}C^{(1)}_{0}\\C^{(2)}_{0}
\end{array}\right)\right\}\\\notag
&\Rightarrow\left(\begin{array}{cccc}p_{11}p_{22} &
-p_{11}p_{12}\\p_{21}p_{22} & -p_{12}p_{21}\end{array}\right)
\left(\begin{array}{cc}X_{m}^{(1)}(t,x) \\
X_{m}^{(2)}(t,x)\end{array}\right)\\\notag
&=\left(\begin{array}{cc}
p_{11}p_{22}X_{m}^{(1)}(t,x)-p_{11}p_{12}X^{(2)}_{m_{2}}(t,x)
\\
p_{21}p_{22}X_{m}^{(1)}(t,x)-p_{12}p_{21}X^{(2)}_{m_{2}}(t,x)\end{array}\right),\
\ t>0,\ x\in \mathbb{R}^{n}.
\hspace{4.5cm}\qed\end{align}
\\
%%%%%%%%%%%%%%% Proof of the Case 4 of Macro system %%%%%%%%%%%%%%%%%%%
                    %%%%%%%%%%%%%%%%%%%%
\noindent {\bf Proofs of Theorems 4 and 5.}\\ The proofs can be
proceeded parallel to the proofs of Theorem 3 and 2, respectively,
and thus we leave them to the reader. \qed\\

%%%%%%%%%%%%%%%%%%%%%%%%%%%%%%%%%%% Reference... %%%%%%%%%%%%%%%%%%%%%%%%%%%%%%%%%%%%%%%%%%%%%%%%%%%%%%%%%%%%%%%%%
\bigskip

\bibliographystyle{plain}
\begin{small}

\end{small}

\end{document}